\documentclass[largeformat]{interact}
\usepackage{epstopdf}
\usepackage[longnamesfirst,sort&compress]{natbib}
\usepackage{algorithm} 
\usepackage{algorithmic} 
\usepackage{multirow} 
\usepackage{xcolor}
\usepackage{setspace}
\usepackage{amsmath}
\usepackage{float}
\usepackage{amsmath}

\theoremstyle{plain}
\newtheorem{theorem}{\textbf{Theorem}}
\newtheorem{lemma}{\textbf{Lemma}}
\newtheorem{corollary}{\textbf{Corollary}}
\newtheorem{remark}{\textbf{Remark}}
\newtheorem{definition}{\textbf{Definition}}
\newtheorem{example}{\textbf{Example}}

\usepackage{url}
\usepackage{subfigure}
\usepackage{fancyhdr}
\usepackage{amsmath}
\usepackage{multirow}
\usepackage{pmat}
\usepackage{arydshln}
\usepackage{color}
\usepackage{graphics} 
\usepackage{graphicx}
\usepackage{indentfirst}

\begin{document}


\title{\LARGE \bf
{{{Input Matrix Construction and Approximation Using a Graphic Approach}}}
}

\author{
\name{Yuan Zhang\textsuperscript{a}\thanks{Emails: zhangyuan14@mails.tsinghua.edu.cn, tzhou@mail.tsinghua.edu.cn}  and Tong Zhou\textsuperscript{a} }
\affil{\textsuperscript{a} Department of Automation, Tsinghua University, Beijing, 100084, P.~R.~China
        }
}

\maketitle

\begin{abstract}
Given a state transition matrix (STM), we reinvestigate the problem of constructing the sparsest input matrix {{with a fixed number of inputs to guarantee controllability of the associated system}}. 
{{A new and simple graph-theoretic characterization is obtained for the sparsity pattern of input matrices to guarantee controllability for a general STM admitting multiple eigenvalues,
 and a deterministic procedure with polynomial time complexity is suggested for constructing real valued input matrices  satisfying the controllability requirement with an arbitrarily prescribed sparsity pattern.}}
  Based on this criterion, some novel results on sparsely controlling a system are obtained.
 It is proven that the minimal number of inputs to guarantee controllability equals the maximum geometric multiplicity of the STM under the constraint that some states can not be directly actuated by inputs (provided such problem is feasible), {{extending the existing results.}} Moreover, the minimal sparsity of an input matrix to ensure system controllability can be affected by the number of independent inputs. Furthermore, a graph-based submodular function is built, leading to a greedy algorithm which efficiently approximates the minimal states needed to be directly actuated by inputs to ensure controllability for general STMs. 
To approximate the sparsest input matrices with a fixed number of inputs, we propose a simple greedy algorithm (non-submodular) and a two-stage algorithm, and demonstrate that the latter algorithm, inspired from graph coloring, has a provable approximation guarantee. 
Finally, we present some numerical results to show the efficiency and effectiveness of our approaches.

\end{abstract}

\section{Introduction}
The synthesis and analysis of large scale systems from a control theoretical perspective have recently been the focus of research with the emergence of complex networks, such as biological transduction networks \citep{nature 2011}, power networks \citep{controllability metrics}, gene regulation networks \citep{Xiongjie_2014}, etc.
Among the related issues, the input/output selection problems are especially important but challenging due to their {{inherent combinatorial nature}} (\citealp{nature 2011, A_Olsehvsky_2014, controllability metrics, T_H_Summers_2016, Tzouma_Control_Effort, zhou_2015,  Y_Zhang_2016, Zhou_minimal_control_2016, Zhou_Automatica}). Researchers have developed various strategies for selecting optimal or suboptimal inputs under the objectives to minimize the number of states that must be directly affected by an actuator{{\footnote{{{A state that is directly affected by an actuator is called an `actuated state' in this paper for simplicity of notion. This state is associated with a nonzero row of the system input matrix.}}}} \citep{A_Olsehvsky_2014, T_H_Summers_2016, Sergio_Pequito_2016},  to have the number of independent inputs as small as possible \citep{nature 2011, Zhou_minimal_control_2016, Zhou_Automatica}, or to optimize certain performance metrics such as the worst-case control energy \citep{controllability metrics}, consensus in the presence of noise \citep{A_Clark_2014}, etc.  See \citet{multi-agent_2012} for surveys.

Particularly, it has been proven that, determining the minimal number of actuated states to ensure controllability is NP-hard \citep{A_Olsehvsky_2014}. Nevertheless, a simple greedy algorithm, which maximizes the rank increase of controllability matrix in each iteration, can achieve a logarithmic approximation factor, which is the best performance obtained in polynomial time. 
By contrast, the minimal number of inputs to guarantee controllability has a closed form solution, which equals the maximum geometric multiplicity of the system state transition matrix (STM)  \citep{Zhou_minimal_control_2016}. All the desirable real valued input matrices with this minimal number of inputs are parameterized explicitly in \citet{Zhou_minimal_control_2016}.  Similar actuator/sensor deployment problems have been investigated in \citet{nature 2011}, \citet{Sergio_Pequito_2016}, etc. under structural framework, i.e., a system whose system matrices are either fixed zeros or free parameters. For example, the minimal number of inputs required for structural controllability is given using the matching theory \citep{nature 2011}, and the problem of determining the minimum actuated states to ensure structural controllability is shown to be solvable in polynomial time \citep{Sergio_Pequito_2016}. 

Apart from {{the aforementioned research}} which focuses on the binary concept of controllability, researchers also develop some energy related metrics to quantify controllability and find proper actuator selection strategies to optimize them \citep{controllability metrics, T_H_Summers_2016, Tzouma_Control_Effort}.  For example, in \citet{T_H_Summers_2016}, submodularity is utilized to pick up actuators to optimize certain controllability Gramian related functions; \citet{Tzouma_Control_Effort} has similar motivations but uses a relaxed control energy metric. These investigations extend the binary concept of controllability to quantitative one, which deepens our insight into network controllability.



These efforts greatly deepen our understanding in how to  control a (network) system incurring as less cost {{(such as the number of independent inputs, actuated states, or the control energy)}} as possible. However, note that given an STM, neither the minimal number of actuated states, nor the minimal number of independent inputs, can provide a complete settlement for the sparsity characterization of an input matrix to meet system controllability, as a state variable can be manipulated by many inputs simultaneously, meanwhile an input can actuate more than one states too. {{We note that the vast majority of the {{literature investigates}} the sparsest input configuration in the case where either a single input can control the whole system (e.g., the system dynamics has no repeated eigenvalues, \citealp{A_Olsehvsky_2014, Sergio_Pequito_2017_robust}), or there is no restriction on the number of inputs (e.g., each input actuates only one state, \citealp{Tzouma_Control_Effort, T_H_Summers_2016}). Although under certain conditions on the system dynamics, these two cases are shown to be mathematically equivalent, little attention is
paid to the general case. In this paper, we reinvestigate the problem of constructing the sparsest input matrix with a fixed number of inputs  to guarantee system controllability {{(i.e., input matrix with a prescribed number of columns)}}, where there is no restriction on the spectra of system dynamics. We further assume that some state variables can not be directly actuated, which is more practical in actual engineering \citep{nature 2011, Zhou_minimal_control_2016}. We both study the problem in fundamental properties and algorithmic perspective.}} 
In particular, we propose two new and efficient algorithms, one for the minimal actuated states selection, and one for the sparsest input matrix with a fixed number of inputs to ensure controllability, both with guaranteed approximation performances.
Our main contributions are as follows. Firstly, we give a new and simple graph-theoretic characterization\footnote{More precisely, some algebraic elements are also involved in this characterization. Since such characterization is presented in terms of the input-state-mode digraph defined in Section 3, we still call it as graph-theoretic in a broad sense.} for the sparsity pattern of input matrices to ensure controllability for general STMs admitting multiple eigenvalues, 
 and provide a deterministic procedure to  construct a real input matrix with an arbitrarily prescribed sparsity pattern to ensure controllability in polynomial time.  With this criterion, the relations among several variants of the minimal controllability problems for general STMs can be easily established. 
  Secondly, we build a graph-based function for general STMs and prove its submodularity, leading to a greedy algorithm to approximate the minimal number of actuated states to render controllability. A prominent property of this algorithm is that it reduces much computation burden compared to the existing controllability Gramian or controllability matrix based algorithms in \citet{A_Olsehvsky_2014} and \citet{T_H_Summers_2016}, while maintaining the same approximation guarantees.  
   Thirdly, we prove that the minimal number of inputs to guarantee controllability equals the maximum geometric multiplicity of the STMs even under the constraint that some states can not be directly actuated, provided that such problem is feasible, extending the existing results \citep{Zhou_minimal_control_2016}. 
 Fourthly, we propose a matroid intersection based greedy algorithm and a two-stage algorithm to approximate the sparsest input configuration with a fixed number of inputs to ensure controllability. The latter algorithm, inspired from graph coloring,  is computationally efficient and has provable worst case performance guarantees.
 By contrast, it is shown by counterexample that the mapping from the additional input links (corresponding to the nonzero entries of an input matrix) to the generic dimension of controllable subspaces is not necessarily submodular. As such, a simple greedy algorithm for the aforementioned problem is not accompanied with performance guarantees (although it often performs well). It is worthwhile to mention that, the above results can be directly extended to the corresponding observability problems by duality between controllability and observability. 

The rest is organized as follows. In Section 2, the { {problem formulation}} and some preliminaries are provided.
Section 3 gives a new graph-theoretic characterization for the sparsity pattern of input matrices to ensure system controllability.
The minimal number of inputs to guarantee controllability with state constraints is discussed in Section 4.
Section 5 provides algorithmic aspects on the sparsest controllability problems from a graph-theoretic perspective, with Section 6 presenting numerical results. The last section gives concluding remarks. 

 {\bf{Notations:}} Let $M$ be a matrix and ${\mathcal{J}}$ be a set of integers, then $M_{{{\mathcal{J}}}}$ (, $[M]_{{{\mathcal{J}}}}$ respectively) represents the submatrix of $M$ comprising the columns (rows, resp.) with indices given by ${{\mathcal{J}}}$. $||M||_0$ denotes the number of nonzero entries of matrix $M$. Denote the unity matrix by $I$, whose dimension is omitted when can be inferred from the context. 
 For a set $S$, $|S|$ is its cardinality. Let $\bar M \in \{0,\ast\}^{m\times n}$ be a structured matrix with dimension $m\times n$, i.e., matrix with entries being either fixed zeros or free parameters, where $\ast$ denote free parameters; moreover, denote by $\bar {\mathcal{M}}$ the set of real matrices with sparsity pattern $\bar M$, i.e., $\bar {\mathcal{M}}=\{P\in \mathbb{R}^{m\times n}: P_{ij}=0 \quad \text{if} \quad \bar M_{ij} =0\}$.  Let superscript $\intercal$ denote the transpose of a matrix. 
 Denote a (directed or undirected) graph  ${\mathcal{G}}$ by ${\mathcal{G}}=({\mathcal{V}},{\mathcal{E}})$ where ${\mathcal{V}}=\{v_1,...,v_n\}$ is the vertex set and ${\mathcal{E}}$ the edge set.  A path in a digraph is a sequence of edges without repeated vertices. Given a collection of disjoint subsets $\{{\mathcal{V}}_1,...,{\mathcal{V}}_k\}$  of ${\mathcal{V}}$, a path $v_1\rightarrow , ..., \rightarrow v_k$ is called a ${\mathcal{V}}_1-,...,-{\mathcal{V}}_k$ path, if $v_1 \in {\mathcal{V}}_1$, ... ,$v_k \in {\mathcal{V}}_k$. Given a graph ${\mathcal{G}}$, let $V({\mathcal{G}})$ be its vertex set and $E({\mathcal{G}})$ its edge set. A $k$-clique is an undirected graph with $k$ vertices and every two distinct vertices of them being adjacent.

\section{{ {Problem Formulation}} and Preliminaries}
Consider the following linear time invariant plant $\bf{\Sigma}$
\begin{equation}
\label{plant Eq}
\dot{x}(t)=Ax(t)+Bu(t),
\end{equation}
where $x(t) \in \mathbb{R}^n$ is the state vector, $u(t) \in \mathbb{R}^l$ is the input vector, $A \in \mathbb{R}^{n\times n}$ and $B \in \mathbb{R}^{n \times l}$ are respectively the STM and input matrix.

In (networked) system designs, an important problem is to find a $B$ to make the system $(A,B)$ controllable, while $B$ meets some sparsity properties. To be specific, consider the following three sparsity objectives respectively:
\begin{itemize}
\item {\bf{{\emph{minimal sparsity controllability problem}} (MSCP)}}: finding the sparsest $B \in \mathbb{R}^{n \times l}$ rendering system controllability, where $l\le n$ is prescribed; 
\item {\bf{{\emph{minimal actuated state controllability problem}} (MACP)}}: finding a $B \in \mathbb{R}^{n \times l}$ with the minimum number of nonzero rows rendering system controllability;
\item {\bf{{\emph{minimal input controllability problem}} (MICP)}}: finding a $B$ with the minimum number of nonzero columns guaranteeing system controllability.
\end{itemize} 
The above three problems are sometimes collectively called {\emph{minimal controllability problems}} in {{references}} \citet{nature 2011}, \citet{A_Olsehvsky_2014}, \citet{Sergio_Pequito_2017_robust}, \citet{Zhou_minimal_control_2016}, etc. Another related problem is to find the sparsest diagonal matrix $B_d\in \mathbb{R}^{n \times n}$ to ensure the controllability of $(A,B_d)$, i.e., {\bf{\emph{the minimal diagonal controllability problem}}}, which also belongs to the MACP with the fixed number of inputs $l=n$. In the following, when referring to the MACP for an STM, the number of inputs is set to be $l=n$ by default,\footnote{In fact, as shown in Theorem \ref{theorem 3}, all $B$ satisfying the MACP with different feasible number of inputs  have the same number of nonzero rows.} if no specific value is assigned to $l$. Besides, we would call $S\subseteq\{1,...,n\}$ an optimal solution to the MACP if $S$ is the set of actuated states of an optimal input matrix for the MACP.  
 \emph{The number of inputs $l$ makes essential difference to the solution configuration of the MSCP} (see Section 5). To distinguish such difference, we call the minimal sparsity controllability problem with a fixed number $l$ of inputs by {\bf{$l$-MSCP}} for abbreviation. 

We also reconsider the MICP when there are `forbidden states', i.e., states that can not be directly affected by inputs. This consideration is due to the physical nature that not all the state variables can  receive the input signals directly. Take the series circuit network shown in Fig. \ref{fig0} as an example, where the physical interpretations of the symbols are given in the bottom of Fig. \ref{fig0}. The goal is to control $i_i(t)$ and $u_i(t)$ for each element of the circuit using as less independent voltage sources as possible.  Define state variables $[i_1(t),u_1(t),...,i_N(t),u_N(t)]^{\intercal}$. Then, the state space model is given as
\[{\tiny{\left[\!\! {\begin{array}{*{20}{c}}
{{\dot i_1}(t)}\\
{{\dot u_1}(t)}\\
{{\dot i_2}(t)}\\
{{\dot u_2}(t)}\\
 \vdots \\
{{\dot i_N}(t)}\\
{{\dot u_N}(t)}
\end{array}}\!\! \right] = \left[\!\! {\begin{array}{*{20}{c}}
{ - {R_1}/{L_1}}&{ - 1/{L_1}}&0&{}&{}&{}&{}&{}\\
{1/{C_1}}&0&{ - 1/{C_1}}&0&{}&{}&{}&{}\\
0&0&{ - {R_2}/{L_2}}&{ - 1/{L_2}}&0&{}&{}&{}\\
0&0&{1/{C_2}}&0&{ - 1/{C_2}}&{}&{}&{}\\
{}&{}&{}&{}&{}& \ddots &{}&{}\\
{}&{}&{}&{}&{}&{}&{{\rm{ - }}{R_N}/{L_N}}&{ - 1/{L_N}}\\
{}&{}&{}&{}&{}&{}&{\rm{0}}&{1/{C_N}}
\end{array}}\!\! \right]\left[ {\begin{array}{*{20}{c}}
{{i_1}(t)}\\
{{u_1}(t)}\\
{{i_2}(t)}\\
{{u_2}(t)}\\
 \vdots \\
{{i_N}(t)}\\
{{u_N}(t)}
\end{array}}\!\! \right] + \left[\!\! {\begin{array}{*{20}{c}}
{{e_1}(t)/{L_1}}\\
0\\
{{e_2}(t)/{L_2}}\\
0\\
 \vdots \\
{{e_N}(t)/{L_N}}\\
0
\end{array}}\!\! \right]}}\] As can be seen, for the MICP associated with the above STM, only state variables $i_i(t)$ can be directly affected by inputs, $i=1,...,N$.
\begin{figure}[H]
\begin{minipage}{1\linewidth}
\centering
\includegraphics[width=4.05in]{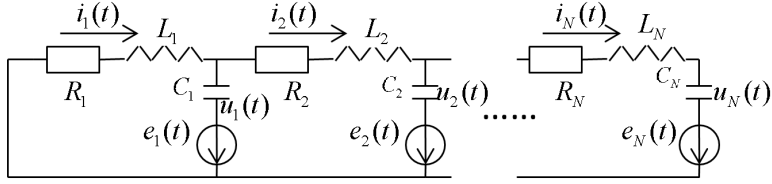}
\begin{center}
\caption{\label{fig0} A series circuit network \citep{Ogata_Modern control_2002}. For each element of the series circuit, $R_i, L_i, C_i$ are the resistance, inductance and the capacitance respectively, $i_i(t),u_i(t)$ are the current of $R_i$ and the voltage of $C_i$ respectively, and $e_i(t)$ is the voltage source which is the possible control input ($i=1,...,N$).} 
\end{center}
\end{minipage}%
\end{figure}

When all eigenvalues of $A$ are distinct, i.e., the system dynamics $A$ is simple, it's shown in \citet{A_Olsehvsky_2014}, \citet{Zhou equivalence}, \citet{Sergio_Pequito_2017_robust}, etc., that the $l$-MSCP is equivalent to the MACP in the sense that they share the same sparsity regardless of the number of inputs ($l\ge 1$). 
 That equivalence essentially lies in the fact that according to the well-known PBH test, the system can be controllable by a vector input matrix with the sparsity pattern $\bar b$ if and only if the support of $\bar b$ (i.e., the set of positions of nonzero entries) intersects with the support of every left eigenvector of $A$ when $A$ is simple \citep{Zhou equivalence, A_Olsehvsky_2014, Sergio_Pequito_2017_robust}. {{However, when $A$ has multiple eigenvalues (i.e., with geometric multiplicities larger than one), each linear combination of those linearly independent left
eigenvectors associated with one eigenvalue is a new eigenvector associated with the same eigenvalue, and the number of distinct supports of eigenvectors can grow exponentially with $n$. 
In such case, the aforementioned condition is only necessary for controllability \citep{Zhou equivalence, Sergio_Pequito_2017_robust}, {{and the complete characterization for sparsity patterns of input matrices to ensure system controllability needs further study}}.  
The main objectives of this paper are as follows: 1) characterizing these minimal controllability problems for a general STM; 2) constructing real valued input matrices with a prescribed sparsity pattern to ensure controllability, especially the MICP in the case where there are forbidden states; 3)  finding efficient algorithms for the MACP and $l$-MSCP for general STMs from a graph-theoretic perspective. 
 Throughout this paper, it's assumed that, an $n\times n$ dimensional real STM $A$ has $p\le n$ distinct eigenvalues, and denote the $i$th eigenvalue by $\lambda_i$, for $i=1,...,p$. Suppose in the $p$ eigenvalues, there are $p_r$ real eigenvalues and $p_c$ complex eigenvalues, which indicates that $p=p_r+p_c$ and $p_c$ is even. 
 Without sacrificing any generality, assume that $\lambda_1, \cdots, \lambda_{p_r}$ are real eigenvalues, and $\{\lambda_{p_r+1}, \lambda_{p_r+1+p_c/2}\},\cdots, \{\lambda_{p_r+p_c/2}, \lambda_{p_r+p_c}\}$ are $p_c/2$ pairs conjugate complex eigenvalues. For $i=1,...,p$, let $k_i$ be the dimension of the left null space of $\lambda_i I- A$; equivalently, ${\rm rank}(\lambda_i I- A)=n-k_i$; that is, $k_i$ is the geometric multiplicity of $\lambda_i$. Denote the maximum geometric multiplicity by $k_{\max}$, i.e., ${k_{\max }} = \mathop {\max }\nolimits_{1 \le i \le p} \left\{ {k_i}\right\}$.  In addition, let $x_{i1},..., x_{ik_i}$ be a set of left eigenvectors of $A$ associated with the eigenvalue $\lambda_i$ which are linearly independent spanning the left null space of $\lambda_i I- A$. Stack these vectors in a matrix $X_i$ as $X_i=[x_{i1},...,x_{ik_i}]$, then $X_i$ is a left eigenbasis of $A$ associated with $\lambda_i$.  Moreover, without losing generality, $\{X_{p_r+i}, X_{p_r+i+p_c/2}\}$ are entry-wise conjugate complex matrices, $i=1,...,p_c/2$. An implicit assumption of this paper is that a collection of eigenbasis of an STM is computationally available.


In the following, we briefly introduce some preliminaries. 

\begin{definition}[Generic rank]\label{definition 1}
 The generic rank of a structured matrix $\bar M$, denoted by ${\rm grank}(\bar M)$,  is the maximum rank it can achieve as the function of its free parameters. 
\end{definition}



The following criterion gives a sufficient and necessary condition for controllability, which is a direct derivation of
the PBH test \citep{multivariable_graph_1988}. 
\begin{lemma}[\citealp{Zhou_minimal_control_2016}]\label{lemma 1}
 Considering the system (\ref{plant Eq}) with $X_i$ being a left eigenbasis of $A$ associated with eigenvalue $\lambda_i$, the system is controllable, if and only if for $i=1,...,p$, $X_i^{\intercal}B$ is of full row rank (FRR).
\end{lemma}


Let $V$ be a finite set and $2^V$ be its power set. A set function $f: 2^{V}\rightarrow {\mathbb{R}}$ assigns a real scalar to each subset of $V$. A nonincreasing function $f: 2^{V}\rightarrow {\mathbb{R}}$ is a set function such that for all $S_1 \subseteq S_2 \subseteq V$, it holds that
$f(S_1)\ge f(S_2)$.

\begin{definition}[Submodularity \citep{Submodular_Wolsey_1982}]\label{definition 2}
  A set function $f: 2^{V}\rightarrow {\mathbb{R}}$ is submodular if for all sets $S_1 \subseteq S_2 \subseteq V$ and any element $s \in V\setminus {S_2}$, it holds that
\begin{equation} \label{submodularity} f(S_1\bigcup \{s\})-f(S_1) \ge f(S_2\bigcup \{s\})-f(S_2). \end{equation}
\end{definition} 



On the basis of nonincreasing function, we have the following criterion to verify a submodular function.

\begin{lemma}[\citealp{Submodular_Wolsey_1982}]\label{lemma 2}
A set function $f: 2^{V}\rightarrow {\mathbb{R}}$ is submodular, if for all $a\in V$, the set function $f_a(S): 2^{V\setminus\{a\}}\rightarrow {\mathbb{R}}$ defined by
$f_a(S)=f(S\bigcup \{a\})-f(S)$
is nonincreasing.
\end{lemma}

Graph coloring is a way of coloring the vertices of an undirected graph such that no two adjacent vertices share the same color.  
\begin{definition}[$k$-coloring, chromatic number \citep{DB_West_graph, graph_coloring}]\label{definition 3}
 A coloring using at most $k$ colors is called $k$-coloring. The smallest number of colors needed to color a graph is called its chromatic number.
\end{definition}

{\emph{Matroids} are combinatorial structures that abstract the notion of linear independence in vector spaces.} A matroid is a pair $(E,\mathcal{I})$ where $E$ is a finite ground set, and
$\mathcal{I}$ is the family of subsets of $E$ which are said to be the independent sets. In this paper, it suffices to think
of a matroid simply as a matrix with respect to the linear independence among its column vectors. For notion of matroids, readers can refer to \citet{Matriod Insertion}, \citet{Matriod Insertion 2}.

\section{A Graph-Theoretic Characterization for Sparsity Pattern of Input Matrices}

In this section, we give a graph-theoretic characterization for the sparsity pattern of $B$ to render system (\ref{plant Eq}) controllable and a deterministic procedure to construct an exact $B$ with given sparsity patterns. To this end, we first define a set of integers $H_i$ as follows:
$${H_i} = \left\{ {{\mathcal{J}} \subseteq \{ 1,...,n\} : X^{\intercal}_{i{{\mathcal{J}}}} \text{ is of FRR},|{\mathcal{J}}| = {k_i}} \right\}$$
for $i=1,...,p$. From this definition, $H_i$ is the collection of 
indices of all $k_i$ linearly independent columns of ${X_i^{\intercal}}$.
\begin{remark}\label{remark 1}
It should be noted that the elements constituting $H_i$ do not vary with the exact $X_i$ that is chosen. A simple algebraic manipulation can interpret this: let ${\tilde X_i}$ be another left eigenbasis of $A$ associated with $\lambda_i$ different from $X_i$. Then, there must exist an invertible matrix $W$ such that
${\tilde X_i} = {X_i}W$. It holds for arbitrary ${\mathcal{J}} \subseteq \{ 1,...,n\}$ that
$$\tilde X^{\intercal}_{i\mathcal{J}}= \tilde X_i^{\intercal}{I_{{\mathcal{J}}}}= {W^{\intercal}}X_i^{\intercal}{I_{{\mathcal{J}}}} = {W^{\intercal}}X^{\intercal}_{i\mathcal{J}}.$$
 As $W^{\intercal}$ is invertible, it can be seen straightforwardly that $\tilde X^{\intercal}_{i{\mathcal{J}}}$ is of FRR if and only if $X^{\intercal}_{i{\mathcal{J}}}$ is of FRR.
\end{remark}

Given an STM $A$, a collection of its left eigenbases $X_i|_{i=1}^p$, and the input sparsity pattern $\bar B \in \{0, \ast\}^{n\times l}$,
define the associated {\emph{input-state-mode (ISM)  digraph}} ${\mathcal{G}}(A,\bar B) = ({\mathcal{V}},{\mathcal{E}})$ as follows. The vertex set ${\mathcal{V}}={\mathcal{U}} \cup {{\mathcal{V}}_s} \cup {{\mathcal{V}}_m}$, where ${\mathcal{U}} = \{ {u_1},...,{u_l}\}$ denotes the input vertices, ${{\mathcal{V}}_s} =\{ 1,...,n\}$ the state vertices, and ${{\mathcal{V}}_m} = \{{m_{{11}}},...,{m_{{1{{k_1}}}}},......,{m_{{p1}}},...,{m_{{p{{k_p}}}}}\}$ the mode vertices respectively. Notice that for each eigenvalue $\lambda_i$, there are $k_i$ mode vertices associated with it, given by $m_i=\{m_{i1},...,m_{i{k_i}}\}$.
The directed edge set  ${\mathcal{E}}={{\mathcal{E}}_{{\mathcal{U}},{{\mathcal{V}}_s}}} \cup {{\mathcal{E}}_{{{\mathcal{V}}_{s}},{{\mathcal{V}}_m}}}$, where ${{\mathcal{E}}_{{\mathcal{U}},{{\mathcal{V}}_s}}} = \{({u_i},{j}):{\bar B_{ji}} \ne 0,{u_i} \in {\mathcal{U}},{j} \in {{\mathcal{V}}_s}\} $ denotes the links from inputs ${\mathcal{U}}$ to states ${\mathcal{V}}_s$, and ${{\mathcal{E}}_{{{\mathcal{V}}_s},{{\mathcal{V}}_m}}} = \{ (i,{m_{{jk}}}):{(X_j^{\intercal})_{ki}} \ne 0, i\in {\mathcal{V}}_s, m_{j_k} \in {\mathcal{V}}_m\}$ the links from states ${\mathcal{V}}_s$ to mode vertices ${\mathcal{V}}_m$. ${(X_j^{\intercal})_{ki}} \ne 0$ means that the $k$th mode associated with $\lambda_j$ can be affected by the inputs injected to  state $i$. {{An illustration of the ISM digraph of the following Example \ref{example 1} is given in Fig. \ref{fig1}.}} 
Notice that the graph representation is noninvariant subject to the exactly chosen $X_i$ (in the structure of ${{\mathcal{E}}_{{{\mathcal{V}}_s},{{\mathcal{V}}_m}}}$), while the results are unaffected by the nonuniqueness of the graph representation due to Remark \ref{remark 1}.   

{{
\begin{example}\label{example 1}
Suppose {{ {\small{$$ A = \left[ {\begin{array}{*{20}{c}}
{4/3}&0&0&{ - 4/3}&0&0\\
0&1&0&0&0&0\\
0&0&3&0&0&0\\
{ - 1/6}&0&0&{5/3}&0&0\\
0&0&{ - 3}&0&2&0\\
0&1&0&0&0&3
\end{array}} \right].$$}}}}The STM $A$ has $p=3$ distinct eigenvalues, which are $\lambda_1=1, \lambda_2=2, \lambda_3=3$ respectively, with geometric multiplicities $k_i=2$ for $i=1,2,3$.
One collection of left eigenbases of $A$ is given by
{{\begin{equation}\label{ism1} \left[ {\begin{array}{*{20}{c}}
{X_1^{\intercal}}\\
{X_2^{\intercal}}\\
{X_3^{\intercal}}
\end{array}} \right]=\begin{pmat}[{.....}]
1&0&0&2&0&0\cr
0&1&0&0&0&0 \cr\-
0&0&3&0&1&0\cr
{-1}&0&0&4&0&0 \cr\-
0&1&0&0&0&2\cr
0&0&1&0&0&0 \cr
\end{pmat}.\end{equation}}}Accordingly, $H_1=\left\{\{1,2\},\{2,4\}\right\}$, $H_2=\left\{ \{1,3\},\{1,5\},\{3,4\},\{4,5\} \right\}$ and $H_3=\left\{\{2,3\},\{3,6\}\right\}$. Considering $\bar B{\rm{ = }} {\footnotesize{ {\left[ {\begin{array}{*{20}{c}}
*&*&0&0&0&0\\
0&*&*&0&0&0
\end{array}} \right]^{\intercal}} }}$, the ISM digraph associated with the eigenbases (\ref{ism1}) is given in Fig. \ref{fig1}.
\end{example}
}}
\begin{figure}
\begin{minipage}[t]{1\linewidth}
\centering
\includegraphics[width=1.6in]{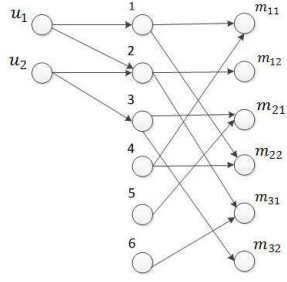}
\begin{center}
\caption{\label{fig1}The input-state-mode digraph ${\mathcal{G}}(A,\bar B) = ({\mathcal{U}}\bigcup {{{\mathcal{V}}_s}} \bigcup {{{\mathcal{V}}_m}}, {{\mathcal{E}}_{{\mathcal{U}},{{\mathcal{V}}_s}}}\bigcup {{{\mathcal{E}}_{{{\mathcal{V}}_s},{{\mathcal{V}}_m}}}})$ of Example \ref{example 1}.}
\end{center}
\end{minipage}%
\end{figure}

\begin{theorem}[characterization of sparsity patterns of input matrices to ensure controllability]\label{theorem 1}
 Given an STM $A \in \mathbb{R}^{n\times n}$ and the sparsity pattern $\bar B \in \{0, \ast\}^{n\times l}$ of input matrices, the following statements are equivalent:
\begin{itemize}
\item[(a)] there exists a real valued $B$ with the sparsity pattern $\bar B$, such that $(A,B)$ is controllable;
\item[(b)] for each $\lambda_i$, $1\le i \le p$, there exists a $k_i\times k_i$ dimensional submatrix of $\bar B$ whose generic rank is $k_i$, and whose row index set belongs to $H_i$;
\item[(c)] for each $\lambda_i$, $1\le i \le p$, in the ISM digraph ${\mathcal{G}}(A,\bar B)$ associated with a collection $X_i|_{i=1}^p$, there exists $k_i$ disjoint ${\mathcal{U}}-{\mathcal{V}}_s-m_i$ paths (denoted by $\mathcal{P}_i$), such that $V(\mathcal{P}_i)\bigcap {\mathcal{V}}_s \in H_i$.
\end{itemize}
\end{theorem}


To prove Theorem \ref{theorem 1}, an intermediate lemma involving the determinant of matrix sums is needed. 
\begin{lemma}[Binet-Cauchy Theorem \citep{R_A_Horn_Matrix}] \label{lemma 4}
Let $M = [{a_1},...,{a_n}]$, $N = {[{b_1},...,{b_n}]^{\intercal}}$, where $a_i, b_i \in \mathbb{C}^m$ for $1\le i \le n$, $m \le n$. Then
\begin{equation}
\label{determiant_0} \det (MN) = \sum\limits_{{\mathcal{J}} \in \Delta } {\det ({M_{\mathcal{J}}}) \det ({{[N]}_{\mathcal{J}}})}\end{equation}
where $\Delta {\rm{ =  }}\{{ {\mathcal{J}}} \subseteq {\rm{ \{ 1,}}...{\rm{,n\} }}: {  |{\mathcal{J}}| = m}\} {\rm{ }}$.
\end{lemma}
{\bf{Proof of Theorem \ref{theorem 1}.}} (a) $\Rightarrow$ (b): 
suppose there exists a $B \in \bar {\mathcal{B}}$, such that $(A,B)$ is controllable. According to Lemma \ref{lemma 1}, this means that $X_i^{\intercal}B$ is of FRR for $i=1,...,p$. Therefore, there exists at least a $k_i\times k_i$  submatrix of $X_i^{\intercal}B$ which is invertible. Denote the column indices of such submatrix by $v_i$ in $X_i^{\intercal}B$. Then according to Lemma \ref{lemma 4},
$$\det ((X_i^{\intercal}B)_{{v_i}}) = \sum\nolimits_{{\mathcal{J}}_i \in \Delta_i }{\det (X_{i{\mathcal{J}_i}}^{\intercal})\det ({[B_{v_i}]_{\mathcal{J}_i}})},$$
where $\Delta_i {{ = \{ {\mathcal{J}}}} \subseteq {{\{ 1,}}...{{,n\} : |{\mathcal{J}}| = }}{{{k}}_i}{{\} }}$.
It can be seen that to make $\det ((X_i^{\intercal}B)_{{v_i}})\ne 0$, at least one ${\mathcal{J}}_i\in \Delta_i$ exists such that ${\det (X_{i{{\mathcal{J}}_i}}^{\intercal}{[B_{v_i}]_{\mathcal{J}_i}})}\ne 0$. It follows that $X_{i{{\mathcal{J}}_i}}^{\intercal}$ and ${[B_{v_i}]_{\mathcal{J}_i}}$ must be invertible simultaneously, which indicates that ${\rm grank}({[B_{v_i}]_{\mathcal{J}_i}})=k_i$, and ${\mathcal{J}}_i\in H_i$ by the definition of $H_i$. 

(b) $\Rightarrow$ (a): suppose that (b) holds. Let $\Gamma  = \{ {\eta _1},...,{\eta _{{k_{nz}}}}\}$ be a set of free parameters in $\bar B$ with $k_{nz} \buildrel \Delta \over = ||\bar B||_0$, and denote the obtained input matrix by $\bar B(\Gamma)$.  For each $\lambda_i$, $1\le i \le p$, let $h_i \in H_i$ be the row indices of a $k_i \times k_i$ submatrix with generic rank $k_i$ in $\bar B$, and ${\phi _i}$ the corresponding column indices, which indicates $|\phi_i|=k_i$. By the condition that $[\bar B_{\phi_i}]_{h_i}$ has full generic rank,
it is not difficult to see that $\det({X_i^{\intercal}\bar B(\Gamma )_{{\phi _i}}})$ is a polynomial of $\{{\eta _1},...,{\eta _{{k_{nz}}}}\}$, and is not identically zero (particularly, let the submatrix obtained by deleting rows indexed by $h_i$ from ${\bar B}(\Gamma)_{\phi_i}$ be zero, then $\det({X_i^{\intercal}\bar B(\Gamma )_{{\phi _i}}})=\det(X^{\intercal}_{h_i})\det([\bar B(\Gamma)_{\phi_i}]_{h_i})\ne 0$).  Consequently, the set of real zeros of $\det({X_i^{\intercal}\bar B(\Gamma )_{{\phi _i}}})$ forms an algebraic variety in $\mathbb{R}^{k_{nz}}$ with zero Lebesgue measure, denoted by ${\cal P}(i,{\phi _i})= \{\Gamma \in {\mathbb{R}^{{k_{nz}}}}: \det({X_i^{\intercal}\bar B(\Gamma )_{{\phi _i}}})=0\}$.  Define a set $\cal{P}$ by
$${\cal P}{\rm{ = }}{\mathbb{R}^{{k_{nz}}}}\backslash \bigcup\nolimits_{i = 1}^p {{\cal P}(i,{\phi _i})}.$$  
As $p$ is countable, $\cal P$ is open and dense in ${\mathbb{R}^{{k_{nz}}}}$. 
That means, there always exists a $\Gamma  = \{ {\eta _1},...,{\eta _{{k_{nz}}}}{\rm{\} }} \in {\cal P}$, such that $X_i^{\intercal}\bar B(\Gamma)$ has a $k_i\times k_i$ invertible submatrix $X_i^{\intercal}\bar B(\Gamma)_{\phi_i}$ for all $i=1,...,p$, indicating that $X_i^{\intercal}\bar B(\Gamma)$ is of FRR. By Lemma \ref{lemma 1}, this leads to the controllability of $(A, \bar B(\Gamma))$.

The equivalence between (b) and (c) is obtained by the following observation. An $n_1\times n_2$ (,$n_1\le n_2$) structured matrix $\bar D$ has full row generic rank, if and only if there are $n_1$ independent free parameters among which any two entries don't locate in the same row and the same column of $\bar D$ \citep{DB_West_graph}. 
Based on this and the construction of the ISM digraph, the statement that there are $k_i$ disjoint ${\mathcal{U}}$-${\mathcal{V}}_s$-$m_i$ paths $\mathcal{P}_i$ while the intersection of ${\mathcal{V}}_s$ and $V(\mathcal{P}_i)$ forms an element of $H_i$, is equivalent to that, there exists a $k_i \times k_i$ submatrix in $\bar B$ with full generic rank while its row indices correspond to a $k_i \times k_i$ invertible submatrix of $X_i^{\intercal}$. Therefore, (c) possesses all the properties hold by (b); and vice versa.  $\hfill\blacksquare$

Considering Condition (c) of Theorem \ref{theorem 1}, it states that $k_i$ disjoint ${\mathcal{U}}-{\mathcal{V}}_s-m_i$ paths must exist, and these paths should pass through certain $h_i\in H_i$, for each $1\le i \le p$, {{which has a similar form to the three-dimensional matching{\footnote{A 3-dimensional matching is a generalization of bipartite matching. Given three disjoint sets $X$, $Y$, $Z$ and a set $T$ which consists of triples $(x,y,z)$ such that $x\in X$, $y \in Y$ and $z \in Z$, a three-dimensional matching is a set $M\subseteq T$ satisfying the following property:  for any two distinct triples $(x_1,y_1,z_1)\in M$ and $(x_2,y_2,z_2)\in M$, we have $x_1\ne x_2$, $y_1\ne y_2$, $z_1 \ne z_2$. }} with an additional constraint that the matched vertices of ${\mathcal{V}}_s$ correspond to $k_i$ linearly independent columns of $X_i^{\intercal}$. With a little abuse of terminology,  we call Condition (c) of Theorem \ref{theorem 1} the \emph{independent-matching condition} for simplicity.}} We say a subset of mode vertices $m^s_i\subseteq m_i$ is independently matched by a set of inputs $\mathcal{U}_s \subseteq \mathcal{U}$, if there exist $|m^s_i|$ disjoint $\mathcal{U}_s-\mathcal{V}_s-m^s_i$ paths, denoted by $\mathcal{P}^s_i$, such that $V(\mathcal{P}^s_i)\cap \mathcal{V}_s \subseteq h_i$ for some $h_i\in H_i$.   It is not difficult to see that, when the STM $A$ is simple, the independent-matching condition collapses to a hitting set{\footnote{Given a collection $\mathcal{C}=\{\mathcal{C}_1,...,\mathcal{C}_m\}$ of subsets of $\mathcal{M}$, a hitting set is a set $\mathcal{C}^{\ast}\subseteq \mathcal{M}$, such that $\mathcal{C}^{\ast}$ intersects with $\mathcal{C}_i$ for $i=1,...,m$.}} (see Section 2). 

{{To show a direct application of Theorem 1, let us revisit Example 1.  For the pair $(A, \bar B)$ in Example 1, it can be found that paths $\{u_1\rightarrow 1\rightarrow m_{11}, u_2\rightarrow 2\rightarrow m_{12}\}$, $\{u_2\rightarrow 3\rightarrow m_{21}, u_1\rightarrow 1\rightarrow m_{22}\}$, $\{u_2\rightarrow 3\rightarrow m_{32},u_1\rightarrow 2\rightarrow m_{31}\}$ satisfy the independent-matching condition for modes associated with $\lambda_1$, $\lambda_2$ and $\lambda_3$ respectively. Therefore, there exists a $B\in \bar {\mathcal{B}}$ making $(A,B)$ controllable; for example, setting all free parameters in $\bar B$ to be $1$ makes a controllable pair $(A,B)$.}}

\begin{remark}[Verifying the independent-matching condition in polynomial time] \label{remark0}
Given the STM $A$, a collection of eigenbases $X_i|_{i=1}^p$ and $\bar B$, the independent-matching condition can be verified by leveraging the \emph{matroid intersection algorithm} \citep{Matriod Insertion} in polynomial time. The matroid intersection problem is to find a largest common independent set in two matroids over the same ground set.  
To be specific, define the ground set $E=\{1,...,n\}$, and the independent sets $\mathcal{I}_i=\{J \subseteq E: {\rm rank}(X^{\intercal}_{iJ})=|J|\}$, and $\mathcal{I}_{\bar B}=\{J \subseteq E: {\rm grank}(\bar B^{\intercal}_{J})=|J|\}$. Then $(E, \mathcal{I}_i)$ and $(E,\mathcal{I}_{\bar B})$ are two {{matroids}}.  It can be verified that the $k_i$ mode vertices associated with $\lambda_i$ are independently matched, if and only if the maximum cardinality of the intersections $\mathcal{I}_i\bigcap\mathcal{I}_{\bar B}$ equals $k_i$.
\end{remark}


Although we have characterized the sparsity pattern of input matrices to meet controllability, it is implicit how to construct such exact input matrices {\emph{in a deterministic way}} (rather than a randomized way). This problem is important not only in system synthesis but also in the computational complexity analysis of general minimal controllability problems \citep{Zhou equivalence}.
In the following, we provide a deterministic procedure (Algorithm \ref{algadd}) with polynomial time complexity to generate an exact input matrix with a given feasible sparsity pattern to ensure controllability.
Algorithm \ref{algadd} is a greedy algorithm overall. The key point is that, in each iteration, to ensure that the modes associated with $\lambda_i$ can be controllable deterministically, we must update at least $k_i$ free parameters of $B$ simultaneously to avoid some hyperplanes which make $X^{\intercal}_iB$ fail to be of FRR, for some $i\in\{1,...,p\}$. Besides, to utilize the fact that a nonzero univariate polynomial with degree $k$ has $k$ zeros, we set all these $k_i$ parameters having the same {{increase}} in each iteration, and choose the {{increase}} that maximizes the number of eigenvalues whose associated mode vertices are totally independently matched (i.e., $m^{\ast}$ in Line 8 of Algorithm \ref{algadd}). 

 \begin{algorithm} 
 {{
\caption{: A deterministic procedure to generate a real input matrix with a given sparsity pattern to ensure controllability} 
\label{algadd} 
\begin{algorithmic}[1] 
\REQUIRE The STM $A$, its left eigenbases $X_i|_{i=1}^p$, the input matrix sparsity pattern $\bar B \in \{0,\ast\}^{n\times l}$ 
\ENSURE A real input matrix $B\in \bar {\mathcal{B}}$ such that $(A, B)$ is controllable 
 \STATE For each $i$, $i\in\{1,...,p\}$, find the $k_i \times k_i$ submatrix in $\bar B$ satisfying Condition (b) of Theorem 1 (if not satisfied, the feasible input matrix does not exist), whose column index set is denoted by $\phi_i$, and determine $k_i$ free parameters in the aforementioned submatrix among which any two don't locate in the same rows and columns, denoting their corresponding unit basis matrices in $\mathbb{R}^{n\times l}$ by $e^{(i)}_1,...,e^{(i)}_{k_i}$ (i.e., $e^{(i)}_j \in \mathbb{R}^{n \times l}$ and the entry in the position of the $j$th free parameter is $1$, $j=1,...,k_i$);
 \STATE   Initialize $B$ to the zero matrix, $c^\ast=0$;
 \FOR {$i=1$ to $p$}
 \IF{$\det(X^{\intercal}_iB_{\phi_i})=0$}
\FOR{$m=1,...,1+\sum\nolimits_{i=1}^p k_i$}
\STATE $B^{(m,i)}=B+\sum\nolimits_{j=1}^{k_i} m e^{(i)}_j$, set $\mathcal{Z}_m=\left\{q: \det(X_q^{\intercal}B^{(m,i)}_{\phi_q})\ne 0, q\in \{1,...,p\}\right\}$;
\ENDFOR
 \STATE Let $m^\ast \in \arg \max_{m} |\mathcal{Z}_m|$, $c^\ast =|\mathcal{Z}_{m^\ast}|$, and $B \leftarrow B+  \sum\nolimits_{j=1}^{k_i} m^\ast e^{(i)}_j$;
\ENDIF
\STATE IF $c^\ast =p$, break and return $B$.
 \ENDFOR
\end{algorithmic}
}}
\end{algorithm}

\begin{theorem}\label{theorem add}
Given $A$ and its associated left eigenbases $X_i|_{i=1}^p$, and  $\bar B$ satisfying the independent-matching condition, Algorithm \ref{algadd} can deterministically find a real valued matrix $B \in \bar {\mathcal{B}}$ such that $(A,B)$ is controllable in polynomial time.
\end{theorem}
{\bf{Proof:}} 
 Denote the free parameters of $\bar B$ corresponding to basis matrices $e^{(i)}_1,...,e^{(i)}_{k_i}$ in Step 1 of Algorithm \ref{algadd} by $\eta^{(i)} _1, ..., \eta^{(i)}_{k_i}$ respectively, for $i=1,...,p$.  Considering the {\bf{for}} loop beginning at Step 3 of Algorithm \ref{algadd}, fix $B$ and $i$ when $\det(X_i^{\intercal}B_{\phi_i})=0$, and let $F_i^{(k)}(m)= \det(X_k^{\intercal}B^{(m,i)}_{\phi_k})- \det(X_k^{\intercal}B_{\phi_k})$ with $B^{(m,i)}$ given in Line 6 of Algorithm \ref{algadd}, $k\in\{1,...,p\}$. From Lemma \ref{lemma 4}, $\det(X_i^{\intercal}B_{\phi_i})= \sum\nolimits_{\mathcal{J}_i\in \Delta_i} \det(X^{\intercal}_{i\mathcal{J}_i})\det([B_{\phi_i}]_{\mathcal{J}_i})$, where $\Delta_i {{ = \{ {\mathcal{J}}: {\mathcal{J}}}} \subseteq {{\{ 1,}}...{{,n\} ,|{\mathcal{J}}| = }}{{{k}}_i}{{\} }}$. 
It can be seen that in $\det(X_i^{\intercal}B_{\phi_i})$, the coefficient of the variable term $\prod\nolimits_{j = 1}^{{k_i}} {\eta _j^{(i)}} $ is nonzero, given by $ \pm \det (X_{i{h_i}}^{\intercal})$, where $h_i\in H_i$ is the row index set of the aforementioned $k_i\times k_i$ submatrix of $\bar B$ in Step 1 of Algorithm \ref{algadd}. Consequently, $F_i^{(i)}(m)$ is a nonzero polynomial of $m$ with degree $k_i$. Therefore, there are at most $k_i$ zeros for $F^{(i)}(m)$, i.e., at most $k_i$ distinct real values exist making $\det(X_i^{\intercal}B^{(m,i)}_{\phi_i})=0$ (noticing $\det(X_i^{\intercal}B_{\phi_i})=0$). 
  Let us consider $q \in \{1,...,p\}\setminus\{i\}$ and $\det(X_q^{\intercal}B_{\phi_q})\ne 0$. Following a similar argument, $F_i^{(q)}(m)$ is either a nonzero univariate polynomial of $m$ with degree at most $k_q$, or identically zero, but can't be some nonzero constant (otherwise it contradicts the fact that $F_i^{(q)}(0)=0$). Therefore, the equation $F_i^{(q)}(m)+\det(X_q^{\intercal}B_{\phi_q})=0$ has at most $k_q$ roots; that is, at most $k_q$ distinct real values can make $\det(X_q^{\intercal}B^{(m,i)}_{\phi_q})=0$. Similarly, for those $r \in \{1,...,p\}\setminus\{i\}$ and $\det(X_{r}^{\intercal}B_{\phi_{r}})= 0$,  $F_i^{({r})}(m)$ is either identically zero, or a univariate polynomial of $m$ with degree at least $1$ and at most $k_{r}$. 

  As a result of the above analysis, for the $i$-th iteration of the $\bf{for}$ loop beginning at Step 3, in any set consisting of $1+\sum\nolimits_{i=1}^p k_i$ distinct real values, there exists at least one real value ${m_0}$, such that: (i) $\det(X_i^{\intercal}B^{({m_0},i)}_{\phi_i})\ne 0$ for fixed $i$ satisfying $\det(X_i^{\intercal}B_{\phi_i})=0$; (ii) $\det(X_q^{\intercal}B^{({m_0},i)}_{\phi_q})\ne 0$ for all $q \in \{1,...,p\}\setminus\{i\}$ satisfying $\det(X_q^{\intercal}B_{\phi_q})\ne 0$; (iii) $\det(X_{r}^{\intercal}B^{({m_0},i)}_{\phi_{r}})\ne 0$ for all $r \in \{1,...,p\}\setminus\{i\}$ satisfying $\det(X_{r}^{\intercal}B_{\phi_{r}})= 0$ while $F_i^{({r})}(m)$ is not identically zero. Since $m^\ast$ is the value maximizing $|\mathcal{Z}_m|$, $m^\ast$ satisfies the above properties simultaneously. Hence, by means of each iteration from Line 4 to Line 9 of Algorithm \ref{algadd}, the number of eigenvalues $\lambda_k$ with $\det(X_k^{\intercal}B_{\phi_k})=0$, $k\in\{1,...,p\}$, reduces at least one; that is, let $\mathcal{Z}_{\hat m^\ast}$ be the corresponding $\mathcal{Z}_{m^\ast}$ in next iteration, then  $|\mathcal{Z}_{\hat m^\ast}|\ge |\mathcal{Z}_{m^\ast}|+1$. After at most $p$ iterations, $|\mathcal{Z}_{m^\ast}|=p$, i.e., $\det(X_i^{\intercal}B_{\phi_i})\ne 0$ for all $i\in \{1,...,p\}$ is satisfied, which, according to Lemma \ref{lemma 1}, certainly leads to the controllability of $(A,B)$, while $B\in {\bar {\mathcal{B}}}$. Step 1 can be implemented in polynomial time using the matroid intersection algorithm, and the rest steps run in polynomial time. Therefore, Algorithm \ref{algadd} has polynomial time complexity.  $\hfill\blacksquare$ 


\begin{remark} From the above proof, the value set in the {\bf{for}} statement in Line 5 of Algorithm \ref{algadd} can be any set consisting of $1+\sum\nolimits_{i=1}^p k_i$ distinct real values, which admits more flexibility in designing an input matrix with entries satisfying certain constraints, such as magnitudes or matrix norm restrictions.  
\citet{Zhou_minimal_control_2016} has provided the parameterization for real input matrices to guarantee controllability. Different from \citet{Zhou_minimal_control_2016}, we set sparsity restriction on the input matrices.  Algorithm \ref{algadd} differs from the procedure in \citet{Zhou equivalence} and the deterministic algorithm in \citet{A_Olsehvsky_2014} in that, we admit the existence of multiple eigenvalues,  which makes the same problem more {{challenging}}. 
 The input matrix construction problem can also been seen as an extension of the so called \emph{matrix completion problem} \citep{Geelen_completion}. 
\end{remark}

In fact, the proof of Theorem 1 makes it clear that controllability is a generic property \citep{Dion_J_M_generic} even for a fixed $A$ and a structured $\bar B$. 
That is, given $A$, almost all numerical instances of $\bar B$ characterized by Theorem \ref{theorem 1} make $(A,B)$ controllable. 
With the availability of Algorithm \ref{algadd}, in the following we focus on the sparsity pattern of input matrices for the associated $l$-MSCP, MACP and MICP rather than an exact numerical instance, which we sometimes call {\emph{input configuration}}.

\section{MICP with Forbidden States}
It has been made clear that, the minimal number of inputs ensuring controllability equals  $k_{\max}$, the maximum geometric multiplicity of the STM, when there is no sparsity pattern restriction on the input matrix,  e.g. \citet{Zhou_minimal_control_2016}.  In the following, as a direct application of Theorem \ref{theorem 1}, we investigate the MICP with the presence of forbidden states. As illustrated in Section 2, it is accepted that this consideration is more practical in actual engineering \citep{Zhou_minimal_control_2016, nature 2011}.

%

To this end, denote the set of  accessible states (i.e., states that can be actuated by inputs) by ${\mathcal{X}}_a\subseteq \{1,...,n\}$. Our main result is given as follows. 

\begin{theorem}\label{theorem 2}
 Given $A$ and an accessible state set ${\mathcal{X}}_a$, the following two statements hold:
\begin{itemize}
\item[(1)] there exist feasible solutions to the MICP with accessible state set ${\mathcal{X}}_a$, only if for each $i=1,...,p$, there exists an $h_i \in H_i$, such that ${h_i} \subseteq {{\mathcal{X}}_a}$;

\item[(2)] if the above condition is satisfied, the minimal number of inputs to ensure system controllability is $k_{\max}$.
\end{itemize}
\end{theorem}
{\bf{Proof.}} From Theorem \ref{theorem 1}, the statement (1) of Theorem \ref{theorem 2} is straightforward.
Now suppose that there is a collection of sets $\{h_1,...,h_p\}$, satisfying $h_i \in H_i$ and $h_i \subseteq {{\mathcal{X}}_a}$ simultaneously for $i=1,...,p$.
Let ${h_i} = \{ {h_{{i_1}}},...,{h_{{i_{{k_i}}}}}\} $, and without losing  generality, assume that ${h_{{i_1}}} < ,..., < {h_{{i_{{k_i}}}}}$. Then, for each $i=1,...,p$,  construct a structured matrix ${{\bar B}^i} \in {\{ 0,\ast \} ^{n \times {k_{\max }}}}$ by letting  ${\bar B}_{{h_{{i_j}}},j}^i = \ast$ for $j = 1,...,{k_i}$, and otherwise zero.
Let $\bar B$ be the entry-wise union of these $\bar B^i$, i.e.,
\[\bar B_{jk} = \bigcup\limits_{i = 1,...,p} {{{\bar B_{jk}}^i}}, j=1,...,n,k=1,...,k_{\max} .\]
It suffices to see that $\bar B$ satisfies the independent-matching condition  and the actuated states are in ${\mathcal{X}}_a$, while $\bar B$ has $k_{\max}$ inputs. 
Hence, by Theorem \ref{theorem 1} and the deterministic procedure Algorithm \ref{algadd}, Statement (2) of Theorem \ref{theorem 2} is proved.  $\hfill\blacksquare$



\begin{example}[Example of the MICP with forbidden states]\label{minimal_input_construct} Considering the circuit in Fig. \ref{fig0} with $N=2$, let the parameters be $R_1=R_2=1$, $L_1=L_2=1$, $C_1=C_2=1$, where their physical units are the standard units and thus omitted. The goal is to control all the state variables $[i_1(t),u_1(t),i_2(t),u_2(t)]$ using the minimal number of independent voltage sources. We have that
{\small{$$A = \left[ {\begin{array}{*{20}{c}}
{{\rm{ - 1}}}&{{\rm{ - 1}}}&{\rm{0}}&{\rm{0}}\\
{\rm{1}}&{\rm{0}}&{{\rm{ - 1}}}&{\rm{0}}\\
{\rm{0}}&{\rm{0}}&{{\rm{ - 1}}}&{{\rm{ - 1}}}\\
{\rm{0}}&{\rm{0}}&{\rm{1}}&{\rm{0}}
\end{array}} \right],$$}}which has a pair of conjugate eigenvalues with $k_1=1, k_2=1$. According to \citet{Zhou_minimal_control_2016}, the minimal number of inputs equals $k_{\max}=1$ and the feasible real input matrices can be parameterized as
$B=T[b_1,b_2,\hat b_1, \hat b_2]^{\intercal}$ where $T$ is the matrix such $J=T^{-1}AT$ is the Jordan canonical form of $A$, $\hat b_1 \ne 0$ (, $\hat b_2$) is the conjugate of $b_1$ (, $b_2$).
However, not all the $[b_1,b_2,\hat b_1, \hat b_2]^{\intercal}$ are feasible in practical, as it may lead to a $B$ having the sparsity pattern as $[*, *, *, *]^{\intercal}$ such that $u_1(t)$ and $u_2(t)$ are directly controlled. Considering the state constraints, let the sparsity pattern of $B$ be $[*, 0, *, 0]$. Then a feasible input matrix is obtained as $B=[0, 0, 1, 0]^{\intercal}$ corresponding to placing the voltage source $e_2(t)$ in Fig. {\ref{fig0}} (for more complicated cases, Algorithm \ref{algadd} can be implemented).
\end{example}

 Another direct application of Theorem \ref{theorem 1} is to show that the minimal number of actuated states to ensure controllability does not vary with the number of inputs for an arbitrary STM $A$ (on the premise that the number of inputs is no less than $k_{\max}$). To this end, similar to the arguments of \citet{A_Olsehvsky_2014} and \citet{Zhou equivalence}, given $A$, we say $A$ is $k$-sparse diagonal controllable, if there exists a diagonal $B_d \in \mathbb{R}^{n\times n}$ with sparsity no more than $k$ such that $(A,B_d)$ is controllable; we say $A$ is $k$-actuated $l$-input controllable, if there is a $B_{l}\in \mathbb{R}^{n\times l}$ whose number of nonzero rows is at most $k$ and number of nonzero columns is at most $l$, such that $(A,B_l)$ is controllable.
\begin{theorem} \label{theorem 3}
Given $A$, $A$ is $k$-sparse diagonal controllable, if and only if $A$ is $k$-actuated $l$-input controllable for any $l \ge k_{\max}$; the same case holds even when there are forbidden states.
\end{theorem}
{\bf{Proof:}} The {\bf{if}} direction is obvious as one can simply set the $i$-th diagonal of $B_d$ to be nonzero if the $i$-th row of $B_l$ is nonzero. The {\bf{only if}} direction follows a similar construction procedure to the proof of Theorem \ref{theorem 2} of constructing an $n\times l$ input matrix from a number of $k$ actuated states, where one just need to replace $k_{\max}$ by $l$ and let $k$ be $\left| {{h_1} \cup ,..., \cup {h_p}} \right|$. The case when there are forbidden states follows a similar argument.  $\hfill\blacksquare$

\section{Approximation Algorithms for MACP and $l$-MSCP: a Graph-Theoretic Approach}
It is proved in \citet{A_Olsehvsky_2014} that the MACP is NP-hard. 
In this section, we develop algorithms to approximate the MACP and $l$-MSCP for general STMs from a graph-theoretic perspective. 
\subsection{Graph-Based Submodular Function for MACP}
In this part we give a new graph-based submodular function for the MACP based on the independent-matching condition, without computing the controllability Gramian or controllability matrix. 

 Let $V=\{1,...,n\}$ be the set of states, and $S\subseteq V$ be the set of states actuated by diagonal inputs. We define a set function $f: 2^V \rightarrow \mathbb{R}$ as the maximum number of mode vertices that can be independently matched by inputs $I_S$ in the associated ISM digraph $\mathcal{G}(A,I_S)$, that is,
\begin{equation}
\label{submodular}f(S) = \sum\nolimits_{i = 1}^p {\mathop {\max}\limits_{{h_i} \in {H_i}} {{ }}} \left| {S\bigcap {{h_i}} } \right|.\end{equation}

It is obvious that, when $A$ is simple, the above function becomes the number of eigenvalues whose associated supports of eigenvectors $X_i$ intersect with $S$.  
According to the definitions of $H_i$ and rank of a matrix, we have $\mathop {\max}\nolimits_{{h_i} \in {H_i}} \left| {S\bigcap {{h_i}} } \right| = {\rm rank}({X}^{\intercal}_{i S})$. Therefore, (\ref{submodular}) can be equivalently rewritten as  
\begin{equation}
\label{submodular2} f(S) = \sum\nolimits_{i = 1}^p {\mathop {\max}\limits_{{h_i} \in {H_i}} {{ }}} \left| {S\bigcap {{h_i}} } \right|=\sum\nolimits_{i = 1}^p {{\rm rank}(X_{iS}^{\intercal})}.\end{equation}

From Theorem \ref{theorem 1}, $f(S)=\sum\nolimits_{i = 1}^p {{k_i}}$ indicates that the resulting system $(A,I_S)$ is controllable. Hence, maximizing $f(S)$ on $S\subseteq V$ leads to controllability. The following theorem reveals the submodularity of $f(S)$.


 \begin{theorem}\label{theorem 4}
 $f(S)$ defined in (\ref{submodular}) is submodular on $S\subseteq V=\{1,...,n\}$. Consequently, the greedy algorithm based on $f(S)$ (Algorithm \ref{alg1}) for the MACP achieves an $\mathcal{O}(log(n))$-approximation, or more precisely, a $log(\sum\nolimits_{i = 1}^p {{k_i}})$-approximation.
 \end{theorem}
{\bf{Proof.}} 
 For all $a \in V$, define a function ${f_a}:{2^{V\backslash \{ a\} }} \to \mathbb{R}$ as
$${f_a}(S) = f({{S}}\bigcup {\{ a\} } {{ }}) - f(S).$$
For each $i\in\{1,...,p\}$, define $f_a^{(i)}: 2^{{V}\setminus \{a\}}\rightarrow \mathbb{R}$ as  $f_a^{(i)}(S) = \mathop {\max }\nolimits_{{h_i} \in {H_i}} \left| {\{ S \cup \{ a\} \} \bigcap {{h_i}} } \right| - \mathop {\max }\nolimits_{{h_i} \in {H_i}} \left| {S\bigcap {{h_i}} } \right|$. Then, it's clear that ${f_a}(S) = \sum\nolimits_{i = 1}^p {f_a^{(i)}(S)}$.
From (\ref{submodular2}), it holds that
\[\begin{array}{l}
f_a^{(i)}(S) = \mathop {\max }\nolimits_{{h_i} \in {H_i}} \left| \{S \cup \{ a\}\} \bigcap {{h_i}}  \right| - \mathop {\max }\nolimits_{{h_i} \in {H_i}} \left| {S\bigcap {{h_i}} } \right|\\
 = {\rm{rank}}(X_{iS \cup \{ a\} }^{\intercal}) - {\rm{rank}}(X_{iS}^{\intercal})\\
 = {\rm{rank}}(X_{iS}^{\intercal}) + {\rm{rank}}(X_{i\{ a\} }^{\intercal}) - {\rm{dim}}({\rm{span}}(X_{iS}^{\intercal})\bigcap {{\rm{span}}(X_{i\{ a\} }^{\intercal})} ) \\
  \quad - {\rm{rank}}(X_{iS}^{\intercal})\\
 = {\rm{rank}}(X_{i\{ a\} }^{\intercal}) - {\rm{dim}}\left({\rm{span}}(X_{iS}^{\intercal})\bigcap {{\rm{span}}(X_{i\{ a\} }^{\intercal})}\right),
\end{array}\]
where $\rm{dim}(\cdot)$ denotes the dimension of a linear space, and $\rm{span}(\cdot)$ the span of column vectors of a matrix. The third equality holds due to the fact that ${\rm{dim}}(y_1\bigcup y_2)={\rm{dim}}(y_1)+{\rm{dim}}(y_2)-{\rm{dim}}(y_1\bigcap y_2)$ for two linear subspaces $y_1$ and $y_2$.
The above relation indicates that $f_a^{(i)}(S)$ is nonincreasing on $S$. 
Hence, ${f_a}(S) = \sum\nolimits_{i = 1}^p {f_a^{(i)}(S)}$ is also nonincreasing with $S$ for any $a \in V$. From Lemma \ref{lemma 2}, $f(S)$ is submodular on $S \subseteq V$.
The remaining statement of Theorem \ref{theorem 4} follows immediately from the submodularity  of $f(S)$ \citep{Submodular_Wolsey_1982}. $\hfill\blacksquare$ 

\begin{algorithm} 
{{
\caption{: Approximation algorithm for the MACP} 
\label{alg1} 
\begin{algorithmic}[1] 
\REQUIRE The STM $A$
\ENSURE   Approximation of the minimal actuated state  set $S$ such that $(A, I_S)$ is controllable 
\STATE   Calculate $X_i$ of $A$, for $i=1,...,p$;
 \STATE   Initialize $V \leftarrow \{1,...,n\}$, $S  \leftarrow \emptyset$;
 \WHILE{$f(S ) < \sum\nolimits_{i = 1}^p {{k_i}} $}
\STATE ${s } \leftarrow a' \in \arg {\kern 1pt} {\kern 1pt} \max {\kern 1pt} {{\kern 1pt} _{a \in {V}\backslash S }}{\kern 1pt} {\kern 1pt} f(S  \cup \{ a\} ) - f(S )$  ($f(S)$ is calculated according to (\ref{submodular2}));
\STATE $S  \leftarrow S  \cup \{ {s }\}$;
\ENDWHILE
\RETURN $S$.
\end{algorithmic}
}}
\end{algorithm}


It is worthwhile to mention that the controllability matrix and the controllability Gramian based submodular functions, which both measure the dimension of controllable subspaces, are used to approximate the MACP, e.g. in \citet{A_Olsehvsky_2014}, \citet{T_H_Summers_2016}. Let us compare Algorithm \ref{alg1} with them in terms of computation complexity.  As is pointed out in \citet{A_Olsehvsky_2014}, the calculation of the controllability matrix $\mathcal{C}(A,I_S)=[I_S, AI_S,...,A^{n-1}I_S]$ is computationally burdensome, and the determination of ${\rm rank}([\mathcal{C}(A,I_S)])$ may encounter numerical instability when the dimension of $A$ becomes large. In each iteration of the greedy algorithm, the controllability Gramian $W_S$, can be determined numerically by solving the {{Lyapunov}} equation
$ \label{Laypunov equaiton} AW_S + W_S{A^{\intercal}} + {I_S}(I_S)^{\intercal} = 0,$ with complexity $\mathcal{O}(n^3)$ when $A\in \mathbb{R}^{n\times n}$ is stable \citep{Spare_Laypunov_2016}. 
In addition, to obtain the rank of $W_S$, it incurs complexity $\mathcal{O}(n^3)$ using the singular value decomposition. While in Algorithm \ref{alg1}, one just needs to calculate the eigenbases of $A^{\intercal}$ for once, which costs $\mathcal{O}(n^3)$ complexity \citep{R_A_Horn_Matrix}. Then, in each iteration, the calculation of ${\rm rank}(X^{\intercal}_{iS})$ can be implemented within linear complexity $\mathcal{O}(n)$, since usually $k_i\ll n$. Hence, computing $f(S)$ incurs at most $\mathcal{O}(pn)$ ($\rightarrow \mathcal{O}(n^2)$) complexity. What's more, once we obtain $X_i|_{i=1}^p$, the determination of $f(S)$ can be implemented parallelly due to the additivity of $f(S)$ (while ${\rm rank}(W_S)$ can't). 
In summary, compared to the controllability {{Gramian}} or the controllability matrix based algorithms,  Algorithm \ref{alg1} incurs much less computation burden. 
\emph{The restriction of  Algorithm \ref{alg1} for a prescribed large-dimensional STM lies in the precision of calculating eigenbases, where round-off errors may influence the results.} See \citet{Sergio_Pequito_2017_robust} for such discussions.   
 Nevertheless, noting that (\ref{submodular2}) builds a bridge between the PBH test and submodular functions,  Algorithm \ref{alg1} can be applied to the input selection for a networked multi-input-multi-output system described in \citet{zhou_2015}, \citet{Stability_Zhou}, \citet{Y_Zhang_2016} without calculating the lumped STM and the associated eigenbases. Details are omitted due to space consideration.

{\begin{remark}[generalization of $f(S)$] Given an STM $A$ and a collection of vector inputs $B=[b_1,...,b_l]$, the idea of Theorem \ref{theorem 4} can also be generalized to select the minimal number of columns $B_S$ from $B$ such that $(A,B_S)$ is controllable. For this problem, $f'(S)=\sum\nolimits_{i = 1}^p {{\rm rank}(X_{i}^{\intercal}B_S)}$ is a feasible objective function. Following Theorem \ref{theorem 4}, it can be validated that $f'(S)$ is submodular on $S\subseteq \{1,...,l\}$.
\end{remark}}
\subsection{Non-equivalence between the $l$-MSCP and MACP}
It is known from \citet{A_Olsehvsky_2014}, \citet{Zhou equivalence}, \citet{Sergio_Pequito_2017_robust} that the optimal solutions to the MACP always have the same sparsity as those to the $l$-MSCP when the STM is simple.  In the following, we show via a simple counterexample that the MACP (or the minimal diagonal controllability problem) is not necessarily equivalent to the $l$-MSCP in the multiple eigenvalue case, \emph{where `equivalent' means that the optimal values of their solutions are equal to each other}.  That means,  the number of inputs might affect the minimal sparsity of an input matrix (\emph{minimal input sparsity}) to ensure controllability.  Then, we explore conditions under which such equivalence holds. This investigation is significant in determining solution to the $l$-MSCP from that of the MACP, as the latter problem is known to be approximated by simple greedy algorithms with the best guarantees in polynomial time \citep{A_Olsehvsky_2014}. Once we have guaranteed the equivalence between the MACP and the $l$-MSCP for a class of STMs $A$, we can construct a solution from the MACP to the $l$-MSCP using some standard manipulations (see Theorem \ref{theorem 6}). 

\begin{example}[Example showing non-equivalence between the $l$-MSCP and MACP]\label{example 3}
  Let us revisit the STM $A$ in Example \ref{example 1}. Let the prescribed number of inputs $l=2$ (satisfying $l\ge k_{\max}=2$). From the ISM digraph in Fig. \ref{fig1}, it is not difficult to see that at least $3$ state variables need to be actuated to ensure controllability; for example, $3$ independent inputs actuating the state set $\{1,2,3\}$ or $\{2,3,4\}$ satisfies the independent-matching condition. However, for $l=2$, it can be validated that any input configuration consisting of $3$ input links with $2$ independent inputs can't satisfy the independent-matching condition. At leat $4$ input links are needed for $l=2$ to ensure controllability; for instance, putting  $1$ on the $(1,1),(2,1),(2,2)$ and $(3,2)$th entries of $\bar B$ is feasible, which is exactly the case shown in the ISM digraph of Fig. \ref{fig1}.
\end{example}


In the following, we leverage graph coloring to characterize the non-equivalence between the $l$-MSCP and MACP. To this end, define an auxiliary graph ${\mathcal{G}}(h_1,...,h_p)=({\mathcal{V}}_p,{\mathcal{E}}_p)$ associated with a given collection $\{h_1,...,h_p\}$ where $h_i\in H_i$ for $i=1,...,p$:  the vertex set ${\mathcal{V}}_p = \bigcup\nolimits_{i = 1,...,p} {{h_i}}$, and the undirected edge set ${\mathcal{E}}_p = \{ (j,k):j,k \in {h_i}, {\text{and}} {\kern 1pt} {\kern 1pt} j \ne k, i = 1,...,p\}$ (no parallel edge is allowed). That is, ${\mathcal{G}}(h_1,...,h_p)$ is the union of $p$ cliques $\mathcal{C}_{h_1}$,...,$\mathcal{C}_{h_p}$, where $\mathcal{C}_{h_i}$ is a $k_i$-clique formed by vertices of $h_i$, $i=1,...,p$. See Fig. \ref{fig_color} for illustration. 

\begin{theorem} \label{theorem 5}
 Given an STM $A$, the associated $l$-MSCP is equivalent to the MACP, if and only if
there exists a collection $\{h_1,...,h_p\}$ where $h_i\in H_i$ for $i=1,...,p$ and ${\bigcup\nolimits_{1 \le i \le p} {{h_i}} }$ is an optimal solution to the MACP, such that the corresponding auxiliary graph ${\mathcal{G}}(h_1,...,h_p)$ has an $l$-coloring.
\end{theorem}
{\bf{Proof.}} Necessity: If such equivalence holds, let $k_{opt}$ be minimal number of actuated states in MACP. Suppose that the sparsity of the optimal value to the $l$-MSCP is also $k_{opt}$ with the corresponding input matrix being $B_l$. Then, there is only one nonzero entry in every nonzero row of $B_l$; otherwise, the number of nonzero rows of $B_l$ is less than $k_{opt}$, which, according to Theorem \ref{theorem 3}, indicates that the system is $k^\ast$-sparse diagonal controllable with $k^{\ast}<k_{opt}$, causing a contradiction. 
Therefore, each state is actuated by at most one input in the optimal solution of the $l$-MSCP.  Let $h_1\in H_1$,...,$h_p \in H_p$ be the collection of actuated states associated with an optimal solution to the $l$-MSCP satisfying the independent-matching condition, with $\left| {\bigcup\nolimits_{1 \le i \le p} {{h_i}} } \right|=k_{opt}$. Then, ${\bigcup\nolimits_{1 \le i \le p} {{h_i}} }$ is also an optimal solution to the MACP. According to the definition of ${\mathcal{G}}(h_1,...,h_p)$, the fact that each state with index in ${\bigcup\nolimits_{1 \le i \le p} {{h_i}} }$ is actuated by at most one input from the $l$ inputs and no two states with indices in $h_i$ share a common input, indicates that ${\mathcal{G}}(h_1,...,h_p)$ has an $l$-coloring. 


Sufficiency:  when $\mathcal{G}(h_1,...,h_p)$ has an $l$-coloring, where $\{h_1,...,h_p\}$ corresponds to an optimal solution to the MACP with $h_i\in H_i, k_{opt} \buildrel \Delta \over =\left| {\bigcup\nolimits_{1 \le i \le p} {{h_i}} } \right|$, no two states with indices in $h_i$ share a common color from the construction of $\mathcal{G}(h_1,...,h_p)$. Let $\{q_1,...,q_{k_{opt}}\} \buildrel \Delta \over  ={\bigcup\nolimits_{1 \le i \le p} {{h_i}} }$, and let the corresponding coloring indices be $\{c_1,...,c_{k_{opt}}\}$, where $1 \le q_i \le n$, $1\le c_i \le l$ for $i=1,...,k_{opt}$. Construct a matrix $\bar B \in \{0, \ast\}^{n \times l}$ by letting $\bar B_{q_i,c_i}=\ast$ for $i=1,...,k_{opt}$. Then, it suffices that $\bar B$ satisfies the independent-matching condition, while every nonzero row of $\bar B$ has only one nonzero entry. Thus, ${\left\| {\bar B} \right\|_0} = {k_{opt}}$,  which, according to Theorem \ref{theorem 3}, indicates that $\bar B$ must be an optimal input configuration of the $l$-MSCP. $\hfill\blacksquare$
 It is probably hard to verify the condition of Theorem \ref{theorem 5} in general, as the optimal solution to the MACP is needed therein, which is NP-hard. Besides, verifying whether a graph has an $l$-coloring is  NP-complete for $l\ge 3$ \citep{DB_West_graph}.  Nevertheless, utilizing properties of graph coloring \citep{graph_coloring}, some easily verified sufficient conditions to guarantee the equivalence between the two problems can be obtained, given as Corollary \ref{corollary 2}.  

\begin{corollary} \label{corollary 2}
Given an STM $A\in \mathbb{R}^{n\times n}$,  under each of the following circumstances, the associated MACP is equivalent to the $l$-MSCP ($l\ge k_{\max}$):  
\begin{itemize}
\item[(i)] every eigenvalue of $A$ has geometric multiplicity $1$, i.e., $k_1=...=k_p=1$;

\item[(ii)] $\sum\nolimits_{i \in \{ 1,...,{p_r} + \frac{{{p_c}}}{2}\} \backslash \{ {i^*}\} } {\scriptsize{\left( \begin{array}{l}
{k_i}\\
2
\end{array} \right)}}  < {k_{\max }}$, where $i^*$ is such that $\lambda_{i^*}$ is the unique eigenvalue with geometric multiplicity being $k_{\max}$ among $\{\lambda_1,...,\lambda_{{p_r} + \frac{{{p_c}}}{2}}\}$;

\item[(iii)] $l \ge \min \{ \sum\nolimits_{i = 1}^{{p_r} + \frac{{{p_c}}}{2}} {{k_i}} ,1+\sum\nolimits_{i = 1}^{{p_r} + \frac{{{p_c}}}{2}} {\scriptsize{\left( \begin{array}{l}
{k_i}\\
2
\end{array} \right)}} \}.$
\end{itemize}
\end{corollary}

%
%
%

{\bf{Proof.}} Case (i) of Corollary \ref{corollary 2} is obvious as the associated auxiliary graph ${\mathcal{G}}(h_1,...,h_p)$ for any collection $h_i|_{i=1}^p$ becomes a set of isolated vertices, which certainly has an $l$-coloring.  For Cases (ii) and (iii), let $\mathcal{G}(h_1,...,h_p)=({\mathcal{V}}_p, {\mathcal{E}_p})$ be the auxiliary graph for any collection $h_i|_{i=1}^p$ satisfying $h_i\in H_i$ and $h_{p_r+i}=h_{p_r+i+p_c/2}$.  Then, Case (ii) is due to the fact that for any graph ${\mathcal{G}}=({\mathcal{V}}, {\mathcal{E}})$ with chromatic number $\chi ({\mathcal{G}})$, it holds that \citep[Chap. 1]{graph_coloring} $\chi ({\mathcal{G}})(\chi ({\mathcal{G}}) - 1)/2 \le \left| {{\mathcal{E}}} \right|$. Note that, with the condition of Case (ii), it follows
\[\left| {{{\mathcal{E}}_p}} \right| \le \frac{{({k_{\max }}{\rm{ - }}1){k_{\max }}}}{2}{\rm{ + }}{\sum\nolimits_{i \in \{ 1,...,{p_r} + \frac{{{p_c}}}{2}\} \backslash \{ {i^*}\} }} {\scriptsize{\left( \begin{array}{l}
{k_i}\\
2
\end{array} \right)}}
 < \frac{{({k_{\max }}{\rm{ - }}1){k_{\max }}}}{2}{\rm{ + }}{k_{\max }} = \frac{{({k_{\max }}{\rm{ + }}1){k_{\max }}}}{2},\]
which means that ${\mathcal{G}}(h_1,...,h_p)$ has a chromatic number at most $k_{\max}\le l$. For Case (iii),  notice that the maximum vertex degree of $\mathcal{G}(h_1,...,h_p)$ is no more than $|{\mathcal{E}}_p|\le \sum\nolimits_{i = 1}^{{p_r} + \frac{{{p_c}}}{2}}  {\scriptsize{\left( \begin{array}{l}
{k_i}\\
2
\end{array} \right)}}
$, and the vertex number satisfies $|{\mathcal{V}}_p| \le \sum\nolimits_{i = 1}^{{p_r} + \frac{{{p_c}}}{2}} {{k_i}}$. Then, Case (iii) is a direct derivation from the fact that the chromatic number of any graph is at most one more than its maximum vertex degree \citep[Chap. 1, Theo. 11]{graph_coloring}, and is no more than its total number of vertices.  $\hfill\blacksquare$

\begin{remark}\label{remark 2}
Case (ii) of Corollary \ref{corollary 2} is suitable for a class of networks having an eigenvalue with geometric multiplicity much higher than the other ones. This property emerges in several real networks, such as scale-free networks with low average degrees, (undirected) ER random networks with low or high connecting probabilities, or Laplacian complete networks, where the eigenvalue $\lambda=0$ usually has very high geometric multiplicity; see \citet{nature 2011}, \citet{NaturePhysics}.  Case (iii) of Corollary \ref{corollary 2} gives an upper bound for $l$ to guarantee the equivalence between the $l$-MSCP and MACP in terms of geometric multiplicities. It can be seen that, Case (i) is a special case of Case (iii).
 \end{remark}

\subsection{Approximation Algorithms for $l$-MSCP}
In the above section, we have shown the non-equivalence between the $l$-MSCP and the MACP in general. Because of the non-equivalence, we can't directly obtain a solution to the $l$-MSCP from the MACP like the simple dynamics case shown in \citet{A_Olsehvsky_2014}, \citet{Sergio_Pequito_2017_robust} and \citet{Zhou equivalence}. We propose two algorithms to approximate the $l$-MSCP. The first algorithm, {\emph{the simple greedy algorithm}}, is a modification of the greedy algorithm for the MACP using the graph-based function defined in this paper. 
The second algorithm, {\emph{the two-stage algorithm}}, is a combination of algorithms for the MACP and dynamic coloring techniques, which comes with guaranteed performance bounds.


\subsubsection{{\bf{Simple greedy algorithm}}}
A natural idea for the $l$-MSCP is to utilize a greedy algorithm like the MACP, i.e., in each iteration selecting an input link to maximize the increase in the dimension of controllable subspaces. The challenge lies in that, given a structured input matrix $\bar B \in \{0,\ast\}^{n\times l}$, it is hard to use numerical methodologies to determine the generic dimension of the controllable subspaces of $(A,\bar B)$. Motivated from the MACP, the objective function in the simple greedy algorithm can alternatively be a modification of the graph-based function (\ref{submodular2}).
To be specific, for a given $A$ and $\bar B \in \{0,\ast\}^{n\times l}$, let $g(\bar B)$ be the {\emph{maximum size of independently matched mode vertices in the associated ISM diagraph ${\mathcal{G}}(A,\bar B)$}}. According to the relation between generic rank of a structured matrix and its associated digraph, $g(\bar B)$ is equivalently written as \citet{Matriod Insertion 2}
\begin{equation}
\label{puregreedy}
g(\bar B) = \sum\nolimits_{i = 1}^p {{\rm grank}(X_i^{\intercal}\bar B)}.
\end{equation}
Based on the above arguments, the simple greedy algorithm has a similar framework to Algorithm \ref{alg1}, i.e., in each iteration choosing a nonzero entry added to $\bar B$ to maximize the increase of $g(\bar B)$, until $g(\bar B)=\sum \nolimits_{i=1}^p {k_i}$.   Following  Remark \ref{remark0}, ${{\rm grank}(X_i^{\intercal}\bar B)}$ in (\ref{puregreedy}) can be computed deterministically using the matroid intersection algorithm in polynomial time.  


\begin{figure}
\begin{minipage}[t]{1\linewidth}
\centering
\includegraphics[width=3.0in]{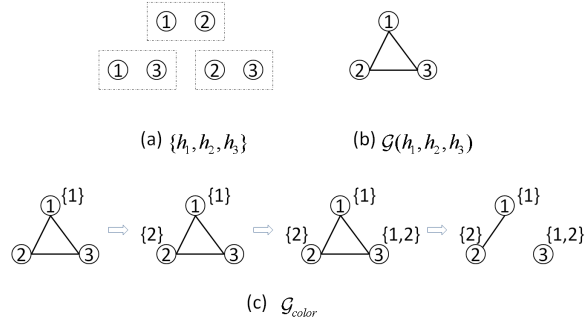}
\begin{center}
\caption{\label{fig_color} (a) and (b): the auxiliary graph $\mathcal{G}(h_1,h_2,h_3)$ with $h_1=\{1,2\}, h_2=\{1, 3\}, h_3=\{2,3\}$ for Example \ref{example 1}. (c): the process of using $2$ colors to color $\mathcal{G}(h_1,h_2,h_3)$, where the numbers in the brackets represent the assigned colors.}
\end{center}
\end{minipage}%
\end{figure}
It should be emphasized that the function $g(\bar B)$ is in general {\emph{not submodular}} on the free parameters of $\bar B$, except for $k_i\equiv 1$.  A counterexample demonstrating the non-submodularity of $g(\bar B)$ is presented in Example \ref{example 4}.  In addition, we note that any functions mapping the additional input links to the generic dimensions of controllable subspaces  are not necessarily submodular. Due to lack of submodularity, it is open to find a nontrivial guaranteed bound for the simple greedy algorithm. In practise, it performs very
well \emph{in terms of approximation} as shown by numerical experiments in Section 6. 

\begin{example}[Non-submodularity of $g(\bar B)$] \label{example 4}
 Consider $A = \left[ {\begin{array}{*{20}{c}}
2&0\\
0&2
\end{array}} \right]$. Let ${\bar B_1} = \left[ {\begin{array}{*{20}{c}}
*&0\\
0&0
\end{array}} \right]$, ${\bar e_{1,2}}{\rm{ = }}\left[ {\begin{array}{*{20}{c}}
0&*\\
0&0
\end{array}} \right]$ and ${\bar B_2} = \left[ {\begin{array}{*{20}{c}}
*&0\\
*&0
\end{array}} \right]$. It can be validated that  
$$ \\~\\~ g({\bar B_1}\bigcup {{\bar e_{1,2}}} ) - g({\bar B_1}) = 0, \\~\\~ g({\bar B_2}\bigcup {{\bar e_{1,2}}} ) - g({\bar B_2}) = 1,$$
consequently, $g({\bar B_2}\bigcup {{\bar e_{1,2}}} ) - g({\bar B_2}) > g({\bar B_1}\bigcup {{\bar e_{1,2}}} ) - g({\bar B_1})$ while ${\bar B_1} \subseteq {\bar B_2}$,
which shows that $g(\bar B)$ is non-submodular.
\end{example}

\subsubsection{{\bf{Two-stage algorithm}}}
To obtain algorithms with guaranteed performances,  we propose the two-stage algorithm (Algorithm \ref{alg3}). This algorithm is motivated by Theorem \ref{theorem 5}, where graph coloring is used to characterize conditions on the equivalence between the $l$-MSCP and the MACP. The basic idea is that, we can obtain a solution to the $l$-MSCP from that of the MACP by adding as few input links as possible to avoid coloring confliction when the number of inputs is limited. However, we strongly suspect that it is  NP-hard to find the smallest difference in sparsity between a feasible solution to the MACP and one to the $l$-MSCP obtained from the MACP.{\footnote{Notice that it is NP-hard to find the smallest number of  colors to color a given graph, under the setting that the number of available colors is fixed and to avoid coloring confliction a vertex can be assigned more than one color simultaneously. If such problem is polynomially solved, then determining whether a graph has a $k$ coloring is polynomially solved too, contradicting to the well-known fact that determining whether graphs admit a $k$ coloring is NP-complete for $k\ge 3$.}}

As suggested by its name, the two-stage algorithm consists of two steps: the first step, using Algorithm \ref{alg1} to approximate the optimal set of actuated states for the associated MACP; the second step, adopting {\emph{dynamic coloring techniques}} \citep{DB_West_graph} to approximate the minimal input links added to avoid coloring confliction caused by limitation of independent inputs.  See
 Fig. \ref{fig_color} for a simple illustration of such process. It should be noted that,  coloring one vertex using more than one colors is admitted in Algorithm \ref{alg3}, which is the key difference from the traditional graph coloring rules. Particularly, the rule (\ref{color_k}) of Algorithm \ref{alg3} is crucial and slightly skillful; see the proof of Theorem \ref{theorem 6}.
 This algorithm has a logarithmic approximation factor if $k_{\max}$ is bounded, and is computationally efficient. Particularly, it is guaranteed that in the cases suggested in Corollary \ref{corollary 2}, Algorithm \ref{alg3} always returns an $\mathcal{O}(log(n))$ approximation. 
\begin{algorithm} 
{{
\caption{: Two-stage algorithm for the $l$-MSCP} 
\label{alg3} 
\begin{algorithmic}[1] 
\REQUIRE The STM $A$, a fixed number of inputs $l$ (, $l \ge {k_{\max }}$) 
\ENSURE   Approximated input configuration for the $l$-MSCP  
\STATE    Calculate $X_i$ of $A$, for $i=1,...,p$;
 \STATE   Find an approximated solution to the associated MACP using Algorithm \ref{alg1}, denoted by $S$, and determine a collection of the corresponding sets $h_i|_{i=1}^{p_r+p_c/2}$, such that $h_i\in H_i$, and $h_i\subseteq S$;
\STATE Construct the auxiliary graph ${\mathcal{G}}(h_1,...,h_{p_r+p_c/2})$;
\STATE  Use $l$ colors to color ${\mathcal{G}}(h_1,...,h_{p_r+p_c/2})$ according to the following rules: index the $l$ colors as $1,2,...,l$, and initialize ${\mathcal{G}}_{color}={\mathcal{G}}(h_1,...,h_{p_r+p_c/2})$. For each iteration, do the following operations, until there is no uncolored vertex in $\mathcal{G}_{color}$:
\begin{itemize}
\item among all uncolored vertices, choose the one which is adjacent to the largest number of differently colored vertices, denoted by $v^\ast$;
\item if vertex $v^{\ast}$ has $l$ differently colored neighbors, assign $k^{\ast}_{\max}$ distinct colors to $v^{\ast}$, where \begin{equation} \label{color_k} k^{\ast}_{\max}=\max \{k_i: v^\ast\in h_i, 1\le i \le p_r+p_c/2\},\end{equation} and remove the edges between $v^{\ast}$ and its neighbors from ${\mathcal{G}}_{color}$; otherwise, assign $v^\ast$ a color different from $v^\ast$'s colored neighbors, such that the number of already used colors is minimized; 
\end{itemize}
\STATE  Map ${\mathcal{G}}_{color}$ to an input configuration $\bar B \in \{0,\ast \}^{n\times l}$: \quad $\bar B_{ij}=\ast$ if vertex $i$ is colored by color $j$ in ${\mathcal{G}}_{color}$, $1\le i \le n$, $1 \le j \le l$; the rest entries of $\bar B$ are fixed zeros.
\end{algorithmic}
}} 
\end{algorithm}

\begin{theorem}[Performance of Algorithm \ref{alg3}]\label{theorem 6}
 The two-stage algorithm  achieves an {\emph{$\mathcal{O}(k_{\max}log(n))$-approximation}} for the $l$-MSCP.  More specifically, let $M^{alg}_{ms}$ be the input configuration for the $l$-MSCP determined by Algorithm \ref{alg3}, $M^{opt}_{ms}$ the optimal input configuration, then one of the following two bounds is guaranteed
\begin{align} \label{firstbound} {\left\| {{M^{alg}_{ms}}} \right\|_0} &\le k_{\max}(\log (N) + 1){\left\| {{M^{opt}_{ms}}} \right\|_0} - (k_{\max}-1)l, \\
 \label{secondbound} {\left\| {{M^{alg}_{ms}}} \right\|_0} &\le (\log (N) + 1){\left\| {{M^{opt}_{ms}}} \right\|_0}\end{align}
where $N = \sum\nolimits_{i = 1}^{{p_r} + \frac{{{p_c}}}{2}} {{k_i}}  \le n$.
In particular, if no vertex is colored by more than one color in the final $\mathcal{G}_{color}$ of Algorithm \ref{alg3}, the bound (\ref{secondbound}) is guaranteed; otherwise the bound  (\ref{firstbound})  is guaranteed. Moreover, under every circumstance of Corollary \ref{corollary 2}, Algorithm \ref{alg3} achieves an $\mathcal{O}(log(n))$-approximation.

\end{theorem}
{\bf{Proof.}}  
We first show the feasibility of the two-stage algorithm. To this end, for the finally obtained ${\cal G}_{color}$ in Step 5 of Algorithm \ref{alg3}, we say a vertex is $k$-colored, if it is assigned with $k$ different colors, $1\le k \le k_{\max}$. Let every vertex of $\mathcal{G}(k_1,...,k_{p_r+p_c/2})$ has the same colors as the corresponding vertex of $\mathcal{G}_{color}$. For each $i\in\{1,...,p_r+p_c/2\}$, considering the subgraph of $\mathcal{G}(h_1,...,h_{p_r+p_c/2})$ induced by $h_i$, denoted by $\mathcal{G}_{h_i}$, its vertices compose of either $1$-colored vertices or $k^{+}$-colored vertices for some $k^+>1$. Recall that a vertex is $k^{+}$-colored only if it has at least $l$ differently $1$-colored neighbors. According to the rule (\ref{color_k}), all those $k^+$ satisfies $k^+ \ge k_i$. As a result, it can be seen that, there exists at least one combination of $k_i$ colors in $\mathcal{G}_{h_i}$, each color chosen from each vertex separately, such that $\mathcal{G}_{h_i}$ is colored with the property that no adjacent vertices share the same color (for example, first choose the unique colors from those $1$-colored vertices, then combinatorially choose one different color from these $k^+$-colored vertices in sequence, $k^+\ge k_i$). That is, for every $i\in\{1,...,p_r+p_c/2\}$, the colored $\mathcal{G}_{h_i}$ admits a $k_i$-coloring. This means that $\bar B$ obtained in Step 5 has the property that
there is a submatrix of $\bar B$ with rows in $h_i$ and having $k_i$ nonzero entries among which every two entries locate in different rows and columns, indicating that such submatrix has generic rank $k_i$. Thus, by Theorem \ref{theorem 1}, the obtained $\bar B$ is a feasible input configuration ensuring controllability. 

For performance bounds, given $A$, denote the optimal diagonal input configuration for the MACP by $M^{opt}_{ma}$, and the approximated one obtained through Algorithm \ref{alg1} (Step 2 of Algorithm \ref{alg3}) by $M^{alg}_{ma}$. By submodularity of the function used in Algorithm \ref{alg1}, it follows
{\begin{equation} \label{AS_IL_relation_1} {\left\| {{M^{alg}_{ma}}} \right\|_0} \le {\left\| {{M^{opt}_{ma}}} \right\|_0(log(N)+1)}  .\end{equation}}Meanwhile, according to Theorem \ref{theorem 3}, the number of nonzero rows of $M^{opt}_{ms}$ is never smaller than that of $M^{opt}_{ma}$, which further leads to
\begin{equation} \label{AS_IL_relation} {\left\| {{M^{opt}_{ma}}} \right\|_0}\le {\left\| {{M^{opt}_{ms}}} \right\|_0}.\end{equation}
Combining (\ref{AS_IL_relation_1}) and (\ref{AS_IL_relation}), it immediately follows that
${\left\| {{M^{alg}_{ma}}} \right\|_0}\le {\left\| {{M^{opt}_{ms}}} \right\|_0}(log(N)+1)$.

Let $k^+$ be some integer larger than $1$.  Note again based on the feasibility analysis in $\mathcal{G}_{color}$,
that a $k^+$-colored vertex emerges only after there exists a vertex who has $l$ differently $1$-colored neighbors. Consequently, if there is no $k^+$-colored vertex in $\mathcal{G}_{color}$, then ${\left\| {{M^{alg}_{ms}}} \right\|_0}={\left\| {{M^{alg}_{ma}}} \right\|_0}\le (log(N)+1){\left\| {{M^{opt}_{ms}}} \right\|_0}$;
if there exist $k^+$-colored vertices, whose number is no more than $\left\| {{M^{alg}_{ma}}} \right\|_0-l$, then
${\left\| {{M^{alg}_{ms}}} \right\|_0}\le k_{\max}({\left\| {{M^{alg}_{ma}}} \right\|_0}-l)+l \le k_{\max}\left({\left\| {{M^{opt}_{ms}}} \right\|_0}(log(N)+1)-l\right)+l,$
which correspond to (\ref{firstbound}) and (\ref{secondbound}) respectively. The above two bounds both guarantee that
${\left\| {{M^{alg}_{ms}}} \right\|_0}\le k_{\max}(log(N)+1){\left\| {{M^{opt}_{ms}}} \right\|_0}.$

The performance bound for Case (i) of Corollary \ref{corollary 2} is obvious. As for Case (ii) of Corollary \ref{corollary 2}, considering the coloring process of Algorithm \ref{alg3}, suppose $k$ colors are already used, $1\le k <l$. Then, the $k+1$th color is needed, only if there is one vertex which has $k$ differently colored neighbors in ${\mathcal{G}}_{color}$; that is, at least $k$ edges exist between such vertex and its colored neighbors. Accordingly, we have that there must exist at least
$\sum\nolimits_{i=1}^{k_{\max}} i={\scriptsize{\left( \begin{array}{l}
{k_{\max}+1}\\
2
\end{array} \right)}}$ edges in ${\mathcal{G}}(h_1,...,h_{p_r+p_c/2})$ if the $k_{\max}+1$th color is used or some $k^+$-colored ($k^+>1$) vertex emerges. However,  in Case (ii), ${\mathcal{G}}(h_1,...,h_{p_r+p_c/2})$ has edges with size at most $\sum\nolimits_{i \in \{ 1,...,{p_r} + \frac{{{p_c}}}{2}\} \backslash \{ {i^*}\} }{\scriptsize{\left( \begin{array}{l}
{k_i}\\
2
\end{array} \right)}} + {\scriptsize{\left( \begin{array}{l}
{k_{\max}}\\
2
\end{array} \right)}}
<{\scriptsize{\left( \begin{array}{l}
{k_{\max}+1}\\
2
\end{array} \right)}}
$, which means that the $k_{\max}+1$th color (if $l\ge k_{\max}+1$) or the $k^+$-colored vertex ($k^+>1$) no longer comes up, indicating that the bound (\ref{secondbound}) is valid. 
The performance bound for Case (iii) of Corollary \ref{corollary 2} follows a similar argument by noting that $\sum\nolimits_{i = 1}^{{p_r} + \frac{{{p_c}}}{2}} {{k_i}} \ge |V({\mathcal{G}}(h_1,...,h_{p_r+p_c/2}))|$ and $\sum\nolimits_{i = 1}^{{p_r} + \frac{{{p_c}}}{2}} {\scriptsize{\left( \begin{array}{l}
{k_i}\\
2
\end{array} \right)}} \ge |E({\mathcal{G}}(h_1,...,h_{p_r+p_c/2}))|$.  $\hfill\blacksquare$

As for complexity of Algorithm \ref{alg3}, Step 1 incurs $\mathcal{O}(n^3)$. In Step 2, the collections $h_i|_{i=1}^{p_r+p_c/2}$ can be determined in a greedy manner: checking the rank increase of $X^{\intercal}_{iS}$ each time when an element is added to $S$, if the gain is $1$ then adding such element to $h_i$. Hence, Step 2 incurs  at most $\mathcal{O}(n^3)$ complexity. In Step 4, the dynamic coloring runs $\mathcal{O}(n)$ iterations and each iteration costs $\mathcal{O}(n^2)$ complexity, which takes $\mathcal{O}(n^3)$ time in total. Steps 3 and 5 can be implemented with linear complexity. To sum up, Algorithm \ref{alg3} has $\mathcal{O}(n^3)$ complexity.

\begin{figure}[H]
  \centering
  \subfigure[]{
    \includegraphics[width=2.37 in]{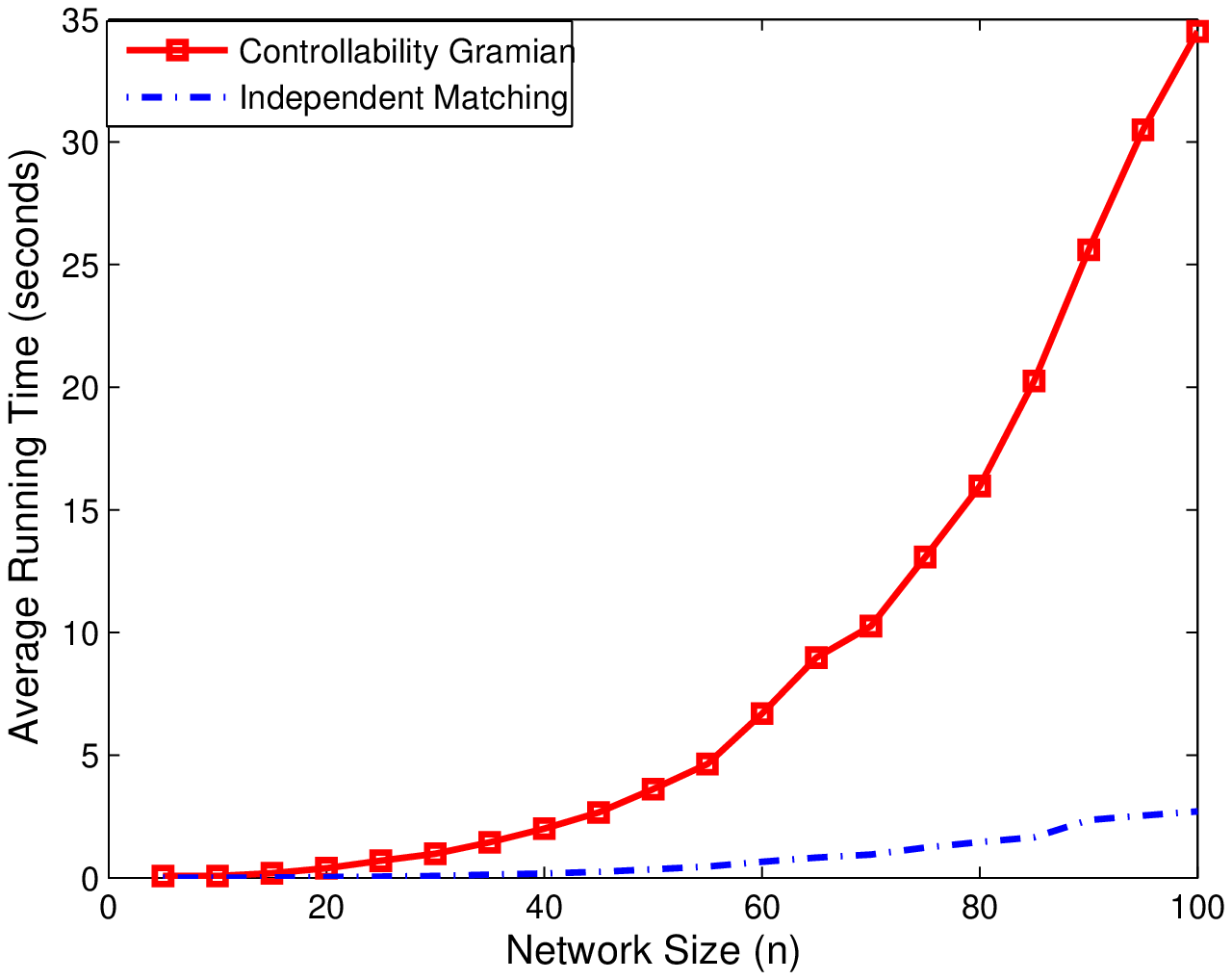}}
  \subfigure[]{
    \includegraphics[width=2.37 in]{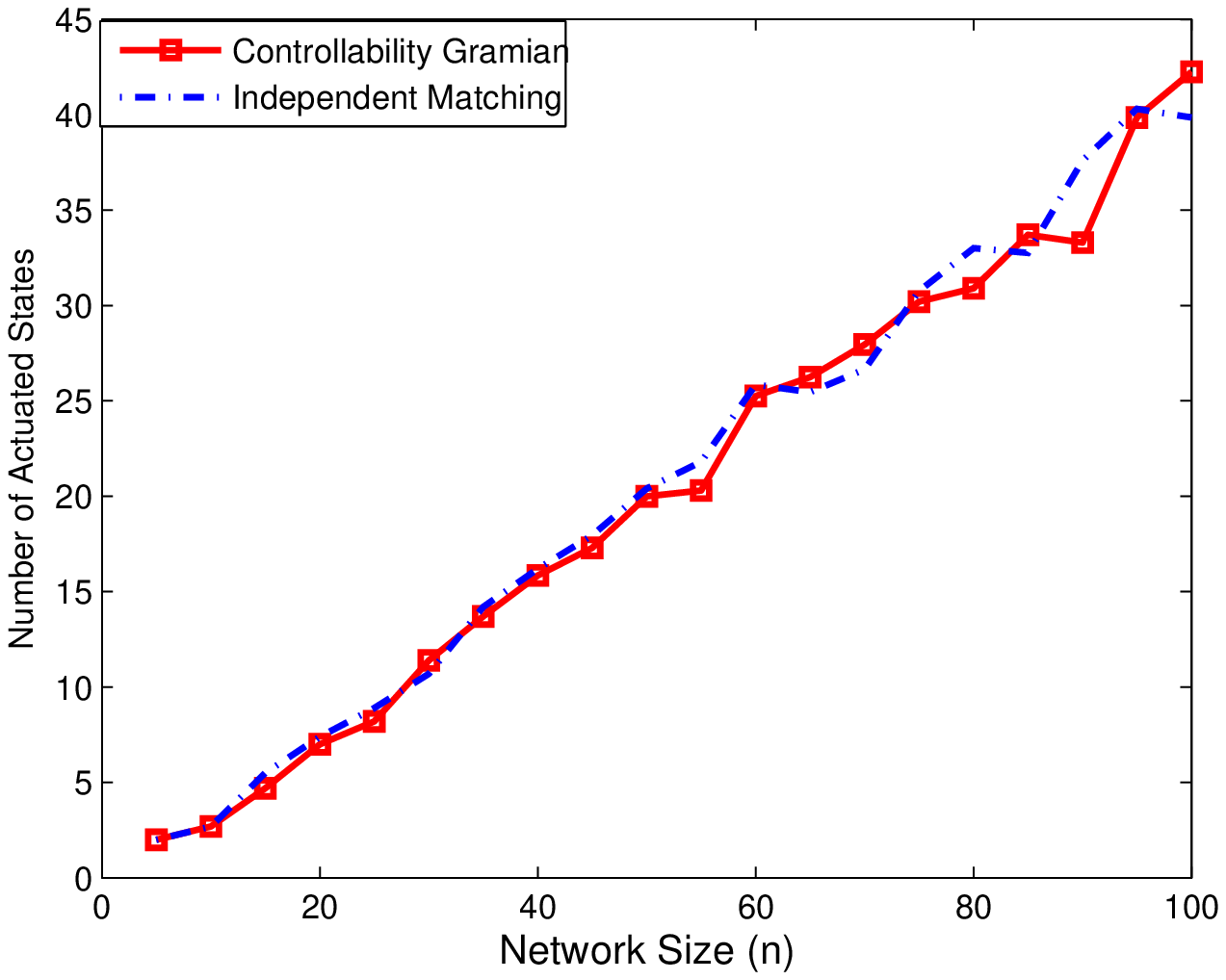}}
  \caption{\label{fig3} Averaged running time (a) and input sparsity (b) of the controllability Gramian based greedy algorithm vs. Algorithm \ref{alg1} for the MACP.}
\end{figure}
\begin{figure}[H]
  \centering
  \subfigure[]{
    \includegraphics[width=2.37 in]{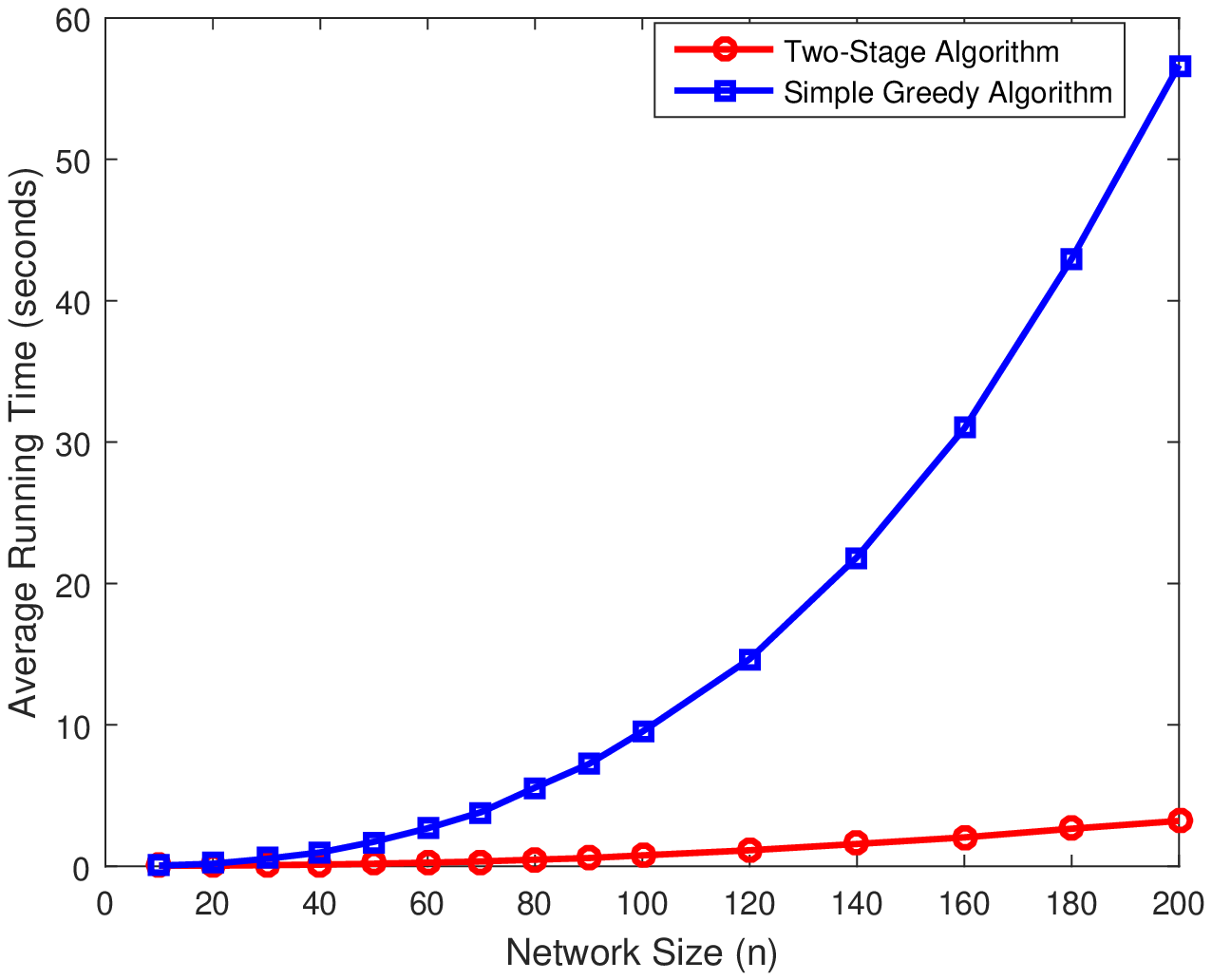}}
  \subfigure[]{
    \includegraphics[width=2.37 in]{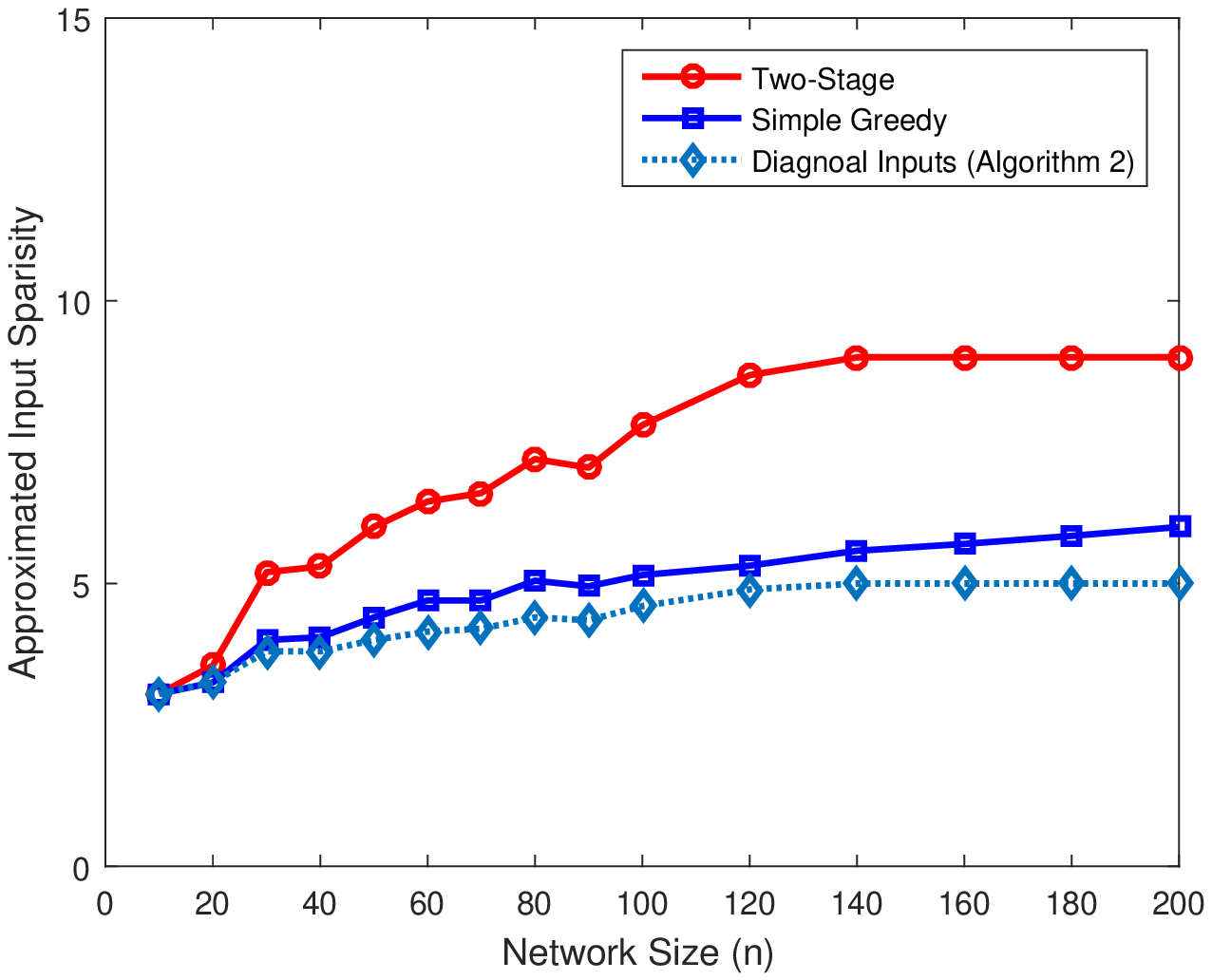}}
  \caption{\label{fig5} Averaged running time (a) and input sparsity  (b) of simple greedy algorithm vs. Algorithm \ref{alg3} for the $l$-MSCP ($l=k_{\max}$).}
\end{figure}

\section{Numerical Simulations}

In the first scenario, we compare the performance between Algorithm \ref{alg1} and the controllability Gramian based greedy algorithm  in \citet{T_H_Summers_2016} for the MACP. For each network size $n$, $n$ ranging from $2$ to $100$, $20$ independent scale free networks are generated via Matlab \citep{ComplexnetworkPackage}, with the power law exponent being $3$ and the average degree proportional to $log(n)$. The weights of the directed edges are uniformly distributed in $[0,1]$. It is worthwhile to mention that, such generated networks tend to have zero eigenvalue with relatively high geometric multiplicity as argued in \citet{nature 2011}, which avoids the trivial case to the largest degree where actuating a single state can control the whole network \citep{A_Olsehvsky_2014}. For each network, the controllability Gramian based greedy algorithm and the independent-matching based algorithm (i.e., Algorithm \ref{alg1}) are used to approximate the corresponding MACP. Each network is stabilized by subtracting $1.1$ times of the largest real part of its eigenvalues for the sake of using {{Lyapunov}} Equation to calculate the controllability Gramian. The average running time and sparsity of the obtained solutions are shown in Figs. \ref{fig3}(a) and (b) respectively. From Fig. \ref{fig3}(a), as the network size grows, the running time of the controllability Gramian based algorithm increases much faster than that of Algorithm \ref{alg1}, while almost the same approximation performances are observed from Fig. \ref{fig3}(b).  These observations are consistent with our theoretical analysis in Section 5.1.  


In the second scenario, we implement simulations to compare the performance between the simple greedy algorithm and Algorithm \ref{alg3} for the $l$-MSCP in non-simple dynamic case. It should be noted that, the benchmark data for real-world autonomous networks which always have eigenvalues with geometric multiplicities higher than one is hard to find. Hence, 
we generate the STMs by inversely using Jordan canonical form decomposition $A=XJX^{-1}$, where $J$ is a block diagonal matrix with Jordan blocks, $X$ is the similarity transformation matrix. In addition, $X$ is set to be invertible{\footnote{To ensure that $X$ is invertible, we transform $X$ to be row diagonally dominant through adding the sum of absolutes of all entries in each row with the addition of $1$ to be the corresponding diagonal entry.}} and sparse with nonzero entries taking up a percentage of $50\%$, and these nonzero entries are randomly located. For each fixed $k_{\max}$ and system size $n$, the geometric multiplicity of the $i$-th distinct eigenvalue $k_i$ is generated  \emph{with equal probability in $\{1,...,k_{\max}\}$  in sequence until  $\sum\nolimits k_i\ge n$}. 
We set $k_{\max}=3$, and for each $n$, $n$ ranging from $20$ to $200$,  generate $20$ independent systems and implement the simple greedy algorithm and Algorithm \ref{alg3} with the number of inputs $l\equiv k_{\max}$. The average running time and input sparsity are shown in Figs. \ref{fig5}(a) and (b) respectively.
It can be seen from Fig. \ref{fig5}(a), as the system size grows, the running time of the simple greedy algorithm increases much faster than that of Algorithm \ref{alg3}. This is due to the fact that the simple greedy algorithm has a much larger search space than Algorithm \ref{alg3}.  In Fig. \ref{fig5}(b), while both algorithms return acceptable approximated input sparsity compared to the system size, the simple greedy algorithms overall achieves a slightly better approximation performance.

Another interesting observation from Fig. \ref{fig5}(b) is that, the approximated minimal input sparsity with $l=k_{\max}$ tends to be larger than that of a diagonal input matrix  for the same $n$. This means that, a possible trade-off may exist between the number of inputs and the minimal input sparsity for systems with non-simple spectra. That is, to ensure system controllability, reducing the number of independent inputs tends to bring the price of increasing the number of input links.
Theorem \ref{theorem 5} indeed theoretically demonstrates the possibility of such trade-off.

\section{Conclusions}
In this paper, we investigate the problem of constructing the sparsest input matrices with a fixed number of inputs to guarantee system controllability. 
Based on the input-state-mode digraph and the matroid intersection, a new and simple graph-theoretic criterion is proposed to characterize the sparsity pattern of input matrices to ensure controllability, along with a deterministic procedure with polynomial time complexity to construct real input matrices with a prescribed sparsity pattern.  From them, a graph-based submodular function is built leading to an efficient greedy algorithm for the MACP.  The MICP is extended to the case where there are forbidden states. Moreover, we show that the number of inputs makes differences to the minimal input sparsity in non-simple dynamic case, i.e., a tradeoff between the minimal input sparsity and the number of inputs may exist, and a simple greedy algorithm to approximate the $l$-MSCP is not accompanied with performance guarantees because of lacking in submodularity. A two-stage algorithm, combining the algorithm for the MACP with techniques in dynamic coloring, achieves acceptable computation efficiency and provable performance guarantees. These results complement and  generalize the results of \citet{A_Olsehvsky_2014}, \citet{Zhou_minimal_control_2016}, etc. 

As further topics, it is interesting to derive nontrivial approximation bounds for the simple greedy algorithm in solving the $l$-MSCP, or to adopt more restrictions on the nonzero entries of input matrices or the STMs, for example variable interdependency e.g. \citet{matrix_net_1982}, \citet{Matriod Insertion 2}, finding computationally efficient algorithms to select inputs and construct controllable matrix pairs.

\section*{Funding}
This work was supported in part by the NNSFC under Grant 61573209 and
61733008.

\end{document}